%version submitted to arXiv

\magnification=1200 \hoffset3cm

\input amstex

\catcode`\@=11\font\myfont=cmr5
\def\myproclaim{%
  \let\savedef@\myproclaim \let\myproclaim\relax
  %\add@missing\endroster \add@missing\enddefinition
  \add@missing\endmyproclaim \envir@stack\endmyproclaim
  \def\myproclaim##1{\restoredef@\myproclaim
    \penaltyandskip@{-100}\preproclaimskip
    {\def\usualspace{\/{\rm\enspace}}%
     \varindent@\rm\ignorespaces##1\unskip
     \frills@{\enspace}}%
    \proclaimfont\ignorespaces}%
  \nofrillscheck\myproclaim}
  \def\endmyproclaim{\revert@envir\endmyproclaim \par\rm
  \vskip\postproclaimskip}
\def\myifdefined#1#2#3{%
  \expandafter\ifx\csname #1\endcsname \relax#2\else#3
    \fi}
\newcount\fornumber\newcount\artnumber\newcount\tnumber
\newcount\secnumber\newcount\subsecnumber
%\newcount\testnumber\newcount\gennumber
%\newcount\stepnumber\newcount\qcount
\newcount\theonumber
\def\Ref#1{%
  \expandafter\ifx\csname mcw#1\endcsname \relax
    \warning{\string\Ref\string{\string#1\string}?}%
    \hbox{$???$}%
  \else \csname mcw#1\endcsname \fi}
\def\Refpage#1{%
  \expandafter\ifx\csname dw#1\endcsname \relax
    \warning{\string\Refpage\string{\string#1\string}?}%
    \hbox{$???$}%
  \else \csname dw#1\endcsname \fi}

\def\warning#1{\immediate\write16{%
            -- warning -- #1}}
\def\CrossWord#1#2#3{%
  \def\x{}%
  \def\y{#2}%
  \ifx \x\y \def\z{#3}\else
            \def\z{#2}\fi
  \expandafter\edef\csname mcw#1\endcsname{\z}
\expandafter\edef\csname dw#1\endcsname{#3}}
\def\Talg#1#2{\begingroup
  \edef\mmhead{\string\CrossWord{#1}{#2}}%
  \def\writeref{\write\refout}%
  \expandafter \expandafter \expandafter
  \writeref\expandafter{\mmhead{\the\pageno}}%
\endgroup}

\openin15=\jobname.ref
\ifeof15 \immediate\write16{No file \jobname.ref}%
\else      \input \jobname.ref \fi \closein15
\newwrite\refout
\openout\refout=\jobname.ref

%defining and labeling sections
\def\dfs#1{\myifdefined{mcwx#1}{}{\warning{multiply defined label #1}}
     \global\advance\secnumber by 1\global\fornumber=0\global\tnumber=0\global\subsecnumber=0
      %\unskip
      \the\secnumber.\if\the\theonumber0
      {\myfont #1}%
\else\relax\fi
     \Talg{#1}{\the\secnumber}
 \expandafter\gdef\csname mcwx#1\endcsname{a}\ignorespaces
}%
\def\dfsubs#1{\myifdefined{mcwx#1}{}{\warning{multiply defined label #1}}
     \global\advance\subsecnumber by 1
      %\unskip
      \the\secnumber.\the\subsecnumber
     \Talg{#1}{\the\secnumber.\the\subsecnumber}
 \expandafter\gdef\csname mcwx#1\endcsname{a}\ignorespaces
}%
%defining and labeling formulas

\def\dff#1{\myifdefined{mcwx#1}{}{\warning{multiply defined label #1}}
     \global\advance\fornumber by 1 \if\the\theonumber0 \text{\myfont #1}%
\else\relax\fi
\unskip (\the\secnumber.\the\fornumber)
\enspace
     \Talg{#1}{(\the\secnumber.\the\fornumber)}
  \expandafter\gdef\csname mcwx#1\endcsname{a} \ignorespaces
}%

%defining and labeling theorems, lemmas, definitions...

\def\dft#1{\myifdefined{mcwx#1}{}{\warning{multiply defined label #1}}%
     \global\advance\tnumber by 1%
       \if\the\theonumber0%
 \ignorespaces    {\myfont #1 }%
\else\relax\fi%
 \the\secnumber.\the\tnumber
    \Talg{#1}{\the\secnumber.\the\tnumber}
  \expandafter\gdef\csname mcwx#1\endcsname{a} \unskip\ignorespaces
}%
\def\dfro#1{\myifdefined{mcwx#1}{}{\warning{multiply defined label #1}}%
     \if\the\theonumber0%
      {\par\noindent\llap{\myfont #1}\par}%
\else\relax\fi%
     \Talg{#1}{(\the\rostercount@)}%
  \expandafter\gdef\csname mcwx#1\endcsname{a} \unskip\ignorespaces
}%
\def\dfr#1{\myifdefined{mcwx#1}{}{\warning{multiply defined label #1}}%
     \if\the\theonumber0%
      {\par\noindent\llap{\myfont #1}\par}%
\else\relax\fi%
     \Talg{#1}{(\the\rostercount@)}%
  \expandafter\gdef\csname mcwx#1\endcsname{a} \unskip\ignorespaces
}%
%defining and labeling articles
\def\dfa#1{\myifdefined{mcwx#1}{}{\warning{multiply defined label #1}}
     \global\advance\artnumber by 1
      \unskip      \the\artnumber%
     \Talg{art.#1}{\the\artnumber}
     \expandafter\gdef\csname mcwx#1\endcsname{a} \unskip\ignorespaces
 }%

\def\rf#1{\ifmmode\Ref{#1}\else $\Ref{#1}$\fi}
\def\df#1{\tag{\myifdefined{mcwx#1}{}{\warning{multiply defined label
#1}}\unskip
     \global\advance\fornumber by 1
     \unskip
    \Talg{#1}{(\the\secnumber.\the\fornumber)}\unskip
 \the\secnumber.\the\fornumber
  \expandafter\gdef\csname mcwx#1\endcsname{a}\ignorespaces
}}%%

\def\donotshowtheoremlabels{\theonumber=1}

%referencing of articles
\def\rfa#1{\Ref{art.#1}}

\csname amsppt.sty\endcsname
\documentstyle{amsppt}

\vsize 18.8 cm

%defining and labeling sections

%\def\ttalg{\tag{\the\secnumber.\the\fornumber}}
%defining and labeling theorems, lemmas, definitions...

\def\dfr#1#2{\myifdefined{mcwx#1}{}{\warning{multiply
defined label #1}}
     \Talg{#1}{#2}
  \expandafter\gdef\csname mcwx#1\endcsname{a} \unskip\ignorespaces
}%

\newdimen\mydim
\long\def\mycenter#1{\mydim=\hsize \advance \mydim by -50pt
\vcenter{\hsize=\mydim
 \noindent\ignorespaces \rm  #1 }}

\def\Li{\roman L}
\def\proof{\demo{Proof}} \def\eproof{\qed\enddemo}
\def\co{\colon}

\def\A{{\bold A}}
\def\N{{\Bbb N}}

\def\R{{\Bbb R}}
\def\D{\roman{d}}
\def\D{\roman{d}}
\def\eps{\varepsilon}
\def\is#1..#2..{\langle #1,#2\rangle}
\def\ci#1..{\mathinner{[{#1}]}}
\def\oi#1..{\mathinner{]{#1}[}}
\def\ro#1..{\mathinner{[{#1}[}}
\def\lo#1..{\mathinner{]{#1}]}}
\def\ici#1..{\mathinner{[\![{#1}]\!]}} \def\Om{\Omega}

\def\Id{\operatorname{\roman{Id}}}
\def\re{\operatorname{\roman{re}}}
\def\Id{\operatorname{\roman{Id}}}
\def\phi{\varphi}
\donotshowtheoremlabels %\donotshowsteps
%\donotshowlabels
%\donotshowformulalabels
\nologo
\hoffset1.5cm

\topmatter
\def\mytime{\the\day.\the\month.\the\year-\the\time}
\title \ Attractors for reaction-diffusion equations on
arbitrary unbounded domains
\endtitle

\author  Martino Prizzi --- Krzysztof P.
Rybakowski
\endauthor
\leftheadtext{\ M. Prizzi --- K. P.
Rybakowski  }
\rightheadtext{\ Reaction-diffusion equations  }

\address Martino Prizzi, Universit\`a degli Studi di
Trieste, Dipartimento di Matematica e Informatica, Via
Valerio, 12, 34127 Trieste, ITALY
\endaddress
\email
prizzi\@ dsm.univ.trieste.it
\endemail

\address Krzysztof P. Rybakowski, Universit\"at Rostock,
Institut f\"ur Mathematik, Universit\"atsplatz 1, 18055
Rostock, GERMANY
\endaddress
\email
krzysztof.rybakowski\@ mathematik.uni-rostock.de
\endemail
%\date March 12, 2006\enddate
\abstract We prove existence of global attractors for parabolic
equations of the form $$\alignedat2  u_t+\beta(x)u-
\sum_{ij}\partial_i(a_{ij}(x)\partial_j u)&=f(x,u),\quad &x\in
\Omega,\,t\in\ro0,\infty..,\\ u(x,t)&=0,\quad &x\in \partial
\Omega,\, t\in\ro0,\infty...
\endalignedat$$ on an arbitrary unbounded domain $\Omega$ in $\R^3$, without
smoothness assumptions on $a_{ij}(\cdot)$ and $\partial\Omega$.
\endabstract

\endtopmatter

\document \hoffset0cm
%\donotshowsteps

%\hsize16cm \vsize25cm
%\baselineskip=14pt
%\donotshowtheoremlabels
\head
\dfs{sec:intr}Introduction
\endhead

In this paper we study the existence of  global
attractors for semilinear parabolic equations of the form $$
 \aligned
  u_t+\beta(x)u-\sum_{ij}\partial_i(a_{ij}(x)\partial_j u) &=f(x,u),\quad x\in
  \Omega,\,t\in\ro0,\infty..,\\ u(x,t)&=0,\quad x\in \partial
  \Omega,\, t\in\ro0,\infty...
 \endaligned\leqno\dff{261105-0951}
$$%
Here, $N=3$ and $\Omega$ is an {\it arbitrary\/} open
set in $\R^N$, bounded or not, $\beta\co\Omega\to\R$ and
$f\co\Omega\times\R\to \R$ are given functions and $\Li u
:=\sum_{ij}\partial_i(a_{ij}(x)\partial_j u)$ is a  linear
second-order differential operator in divergence form. We do not
make any smoothness assumption on $\partial\Omega$ and
$a_{ij}(\cdot)$.

 Notice that, without smoothness assumptions on
$\partial\Omega$ and $a_{ij}(\cdot)$, it is not possible to study
\rf{261105-0951} in the $L^q$ setting for $q\not=2$. The reason
is that one cannot use the regularity theory of elliptic
partial differential equations to characterize the fractional
power spaces generated by $-\Li+\beta(x)$. On the other hand, in order to work  in the $L^2$
setting one must impose growth conditions on $f$. In
particular, for $N=3$ the critical exponent is $\overline\rho=5$. The lack
of regularity also prevents us from being able to use the $\eps$-regular
mild solutions introduced by Arrieta and Carvalho in
\cite{\rfa{Aca}} to treat the critical case. Therefore we shall
assume in this paper  that $f$ has subcritical growth.

There is vast literature concerning  existence of global
attractors for reaction-diffusion equations on bounded domains
(see e.g. \cite{\rfa{Ha},\rfa{La},\rfa{BV1},\rfa{T},\rfa{DCh}}).
In this case the asymptotic compactness property for the solutions
of the equations follows from the compactness of the Sobolew
embedding  $H^1 \subset L^2$.  For unbounded domains this
embedding  is no longer compact and so  new ideas are needed to
obtain the asymptotic compactness property.

In~\cite{\rfa{BV2}} Babin and Vishik  considered an equation of
the form $$
  u_t+u-\Delta u=f(u)+g(x),\quad x\in
  \R^N,\,t\in\ro0,\infty..,
$$ with $f$ satisfying the dissipativeness condition $f(u)u\leq0$
and the monotonicity condition $f'(u)\leq\ell$ . They overcame the
difficulties arising from the lack of compactness by introducing
weighted Sobolev spaces. More recently, Wang considered the same
equation in \cite{\rfa{W}} and established the asymptotic
compactness of the solutions in the space $L^2$, under the same
hypotheses as those in~\cite{\rfa{BV2}}. To this end, he developed a technique based on
tail-estimates of the solutions outside large
balls. The simple remark in~\cite{\rfa{P}} shows that  the
solutions are actually asymptotically compact in the natural energy space
$H^1$.

The equation studied in~\cite{\rfa{BV2},\rfa{W}}  has a very special form. In \cite{\rfa{ACDR}} Arrieta et al. considered the more
general equation $$
 \aligned
  u_t-\Delta u &=f(x,u),\quad x\in
  \Omega,\,t\in\ro0,\infty..,\\ u(x,t)&=0,\quad x\in \partial
  \Omega,\, t\in\ro0,\infty...
 \endaligned
$$ In that paper $\Omega$ is an unbounded domain with
uniformly $C^2$-boundary. The function $f$ has the form
$f(x,u)=m(x)+f_0(x,u)+g(x)$ and satisfies the dissipativeness
condition $f(x,u)u\leq C(x)|u|^2+D(x)|u|$, where $C$ is such that
the semigroup generated by $\Delta +C(x)$ decays exponentially.
The operator $\Delta$ could be replaced by a general second order
differential operator in divergence form like $\Li$, provided the
coefficients $a_{ij}$ are sufficiently smooth. The authors proved
several results about existence of attractors in various Sobolev
spaces, depending on the growth of $f$ and on the summability
properties of $m$, $g$, $C$ and $D$. Their technique is based on
the abstract comparison results of \cite{\rfa{ACaRo}} and
ultimately on the maximum principle for the heat equation.  In
order to apply the comparison results of \cite{\rfa{ACaRo}}, one
needs to check that the nonlinear function $f_0$ satisfies the
following property: for every $r>0$ there exists a constant
$k$ such that the mapping $u(\cdot)\mapsto
f_0(\cdot,u(\cdot))+k u(\cdot)$ is increasing on the ball of
radius $r$ in the functional space in which the problem is set. In
general this property is not satisfied, so one needs to `prepare'
the function $f_0$ before applying the comparison theorem. This
means that one must first find some (local) $L^\infty$-bound for
the solutions and then modify $f_0$ so as to obtain a globally
Lipschitzian function. Such $L^\infty$-bounds are obtained through
a bootstrapping argument which is possible only if
$\partial\Omega$ and $a_{ij}$ satisfy suitable smoothness
assumptions.

In this paper we prove existence of global attractors for the
parabolic equation \rf{261105-0951} on an arbitrary unbounded
domain $\Omega$ in $\R^3$, without smoothness assumptions on
$a_{ij}(\cdot)$ and $\partial\Omega$. To this end we exploit the
tail-estimate technique of Wang and the remarkable fact that the
equation admits a natural Lyapunov functional. Our hypotheses on
the function $f$ are very general and, in particular, they cover the cases considered in \cite{\rfa{ACaRo}}. Moreover, since our proof
does not depend on the maximum principle, it works also for systems
of equations with gradient nonlinearities.

In order to present our results in more detail, let us first describe the notation used in this paper.
\subhead Notation\endsubhead Let  $\Omega$ be an arbitrary open
set in $\R^N$. Given any measurable function $v\co \Omega\to \R$
and any $\nu\in\ro1,\infty..$ we set, as usual, $$
 |v|_{L^\nu}=|v|_{L^\nu(\Omega)}:=\left(\int_\Omega|v(x)|^\nu\,\D
 x\right)^{1/\nu}\le\infty.
$$ Moreover, for $v\in H^1_0(\Omega)$ we set
$|v|_{H^1}=|v|_{H^1(\Omega)}:=(|\nabla
u|_{L^2}^2+|u|_{L^2}^2)^{1/2}$.

We also use the common notation ${\Cal D}(\Omega)$ resp. ${\Cal D}
'(\Omega)$ to denote the space of all test functions on $\Omega$,
resp. all distributions on $\Omega$. If $w\in{\Cal D}'(\Omega)$
and $\varphi\in {\Cal D}(\Omega)$, then we use the usual
functional notation $w(\varphi)$ to denote the value of $w$ at
$\varphi$.

Given  a function $g\co \Omega\times \R\to \R$, we denote by
$\widehat g$  the ({\it Nemitski\/}) operator which associates
with every  function $u\co \Omega\to \R$ the function $\widehat
g(u)\co \Omega\to \R$ defined by
 $$
  \widehat g(u)(x)= g(x,u(x)),\quad x\in \Omega.
 $$
 If $X$ is a normed space and $u\co I\subset \R\to X$ is differentiable into $X$ at $t\in I$ then  we often denote the derivative of $u$ at $t$ by $\partial(u;X)(t)$, in order to indicate its dependence on $X$.

Unless specified otherwise, all linear spaces considered in this
paper are over the real numbers.

\proclaim{Definition \dft{040207-1113}}Let $w\co\Omega\to\R$ be a
measurable function and let $\gamma\in]0,1]$ be a real number. We
say that $w\in{\Cal E}_\gamma$ if and only if one of the following
conditions is satisfied:
 \roster
\item $\gamma\in ]0,1[$ and there exists a constant $C>0$ such
that for all $\eps>0$ and all $u\in H^1_0(\Omega)$
$$\bigl||w|^{1/2} u\bigr|_{L^2}\leq C\bigl(\gamma \eps
|u|_{H^1}+(1-\gamma)\eps^{-\gamma/(1-\gamma)}|u|_{L^2}\bigr)
$$
\item $\gamma=1$ and there exists $C>0$ such that for all $u\in H^1_0(\Omega)$
$$\bigl||w|^{1/2} u\bigr|_{L^2}\leq C |u|_{H^1}.$$
\endroster
\endproclaim

\remark{Remark}Denote by $L^\nu_{\roman u}(\R^N)$  the set of
measurable functions $v\co \R^N\to \R$ such that $$
 |v|_{L^\nu_{\roman u}}:=\sup_{y\in \R^N}\left(\int_{
 B(y)}|v(x)|^\nu\,\D x\right)^{1/\nu}<\infty,
$$ where, for $y\in\R^N$, $B(y)$ is the open unit cube in $\R^N$
centered at $y$. By carefully checking the proof of Lemma 3.2
in~\cite{\rfa{PR1}}, we obtain that if the trivial extension
$\widetilde w$ of $w$ to $\R^N$ lies in $L^\nu_{\roman u}(\R^N)$
for some $\nu\in[ (N/2),\infty[$ then $w\in{\Cal E}_\gamma$, with
$\gamma=(6\nu'-2^*)/4\nu'$. Notice that if $\nu=N/2$ then
$\gamma=1$.\endremark

We make the following assumptions:
\proclaim{Hypothesis~\dft{221105-1756}}\roster
\item\dfro{171205-1239} $a_0$,
$a_1\in\oi0,\infty..$ are constants and $a_{ij}\co
\Omega\to \R$, $i$, $j=1,\dots,N$ are functions in
$L^\infty(\Omega)$ such that $a_{ij}=a_{ji}$, $i$,
$j=1,\dots,N$, and for every $\xi\in\R^N$ and a.e.
$x\in\Omega$, $a_0|\xi|^2\le \sum_{i,j=1}^N
a_{ij}(x)\xi_i\xi_j\le a_1|\xi|^2 $.
$A(x):=(a_{ij}(x))_{i,j=1}^N$,
$x\in\Omega$.\item\dfro{171205-1240} $\beta\co
\Omega\to \R$ is a measurable function with the property
that
  for every
$ \eps\in \oi0,\infty..$ there is a $C_{\eps}\in \ro0,\infty..$
with $\bigl||\beta|^{1/2} u\bigr|_{L^2}^2\le  \eps
|u|_{H^1}^2+C_{\eps}|u|_{L^2}^2$  for all $u\in H^1_0(\Omega)$
\rom(this is slightly less restrictive than the requirement that $\beta\in{\Cal
E}_\gamma$ for some $\gamma\in]0,1[$ \rom).\item
$\lambda_1:=\inf\bigl\{\,\displaystyle\int_\Omega[\sum_{i,j=1}^N
a_{ij}\partial_i u\partial _j u+\beta |u|^2]\,\D x  \bigm| u\in
H^1_0(\Om),\,|u|_{L^2}=1\,\bigr\}>0$.\endroster\endproclaim
\proclaim{Hypothesis~\dft{231105-0903}} \roster\item
 $\overline C$ and $\overline\rho\in\ro0,\infty..$
are constants with $\overline\rho\in\ro2,4..$ and $c\co\Om\to \ro0,\infty..$ is a function
with $c\in L^1(\Om)$;
\item\dfro{020207-1852} $a\co
\Omega\to \R$ is a measurable function with the property that
$a\in{\Cal E}_\gamma$ for some $\gamma\in]0,1[$.
 \item $f\co \Omega\times
 \R\to \R$ is such that $x\mapsto f(x,u)$ is Lebesgue measurable for all $u\in\R$ and $u\mapsto f(x,u)$ is a $C^1$-function for a.e. $x\in \Omega$;
\item $|\partial_u f(x,u)|\le  \overline
C(a(x)+|u|^{\overline\rho})$ for every $u\in\R$ and a.e. $x\in\Omega$;\item $f(\cdot,0)\in L^2(\Omega)$; \item $f(x,u)u\le
c(x)$ and $
\int_0^u f(x,s)\,\D s\le c(x)$ for a.e.  $x\in \Om$ and every
$u\in\R$. \endroster  \endproclaim

Under Hypothesis~\rf{221105-1756} the differential operator $u\mapsto -Lu+\beta(x)u$ defines a positive self-adjoint operator $\A\co D(\A)\subset \bold X\to \bold X$ on the Hilbert space ${\bold X}=L^2(\Omega)$. $D(\A)$ endowed with the graph norm of $\A$ is continuously included in $H^1_0(\Omega)$.
The operator $\A$ generates the family $X^\alpha=D(\A^{\alpha})$, $\alpha\in[0,\infty[$, of fractional power spaces.
Setting $X^{-\alpha}=(X^\alpha)'$, $\alpha\in[0,\infty[$, we can construct a family $\A_{(\alpha)}$, $\alpha\in\R$, of self-adjoint operators, such that $\A_{(\alpha)}\co X^\alpha\to X^{\alpha-1}$. Moreover, $D(\A_{(\alpha)}^{\beta})=X^{\alpha+\beta-1}$ for all $\alpha$, $\beta\in\R$.

 Under Hypothesis~\rf{231105-0903} one can find an $\alpha\in\oi0,1..$ such that the function $f$ generates  a locally Lipschitzian Nemitski operator $\bold f\co H^1_0(\Omega)=D(\A_{(-\alpha+1)}^{\alpha+1/2})\to X^{-\alpha}$.
 By general results on abstract parabolic equations (see e.g. \cite{\rfa{He}}),~\rf{261105-0951} generates a local semiflow $\pi$ on $H^1_0(\Omega)$. The choice of $\alpha$  depends on $\overline\rho$ and $\gamma$.
The semiflow $\pi$ does not depend on the choice of $\alpha$.

The main result of this paper can now be stated as follows.
 \proclaim{Theorem~\dft{010207-0610}} Assume Hypotheses \rf{221105-1756} and \rf{231105-0903}. Then
 $\pi$ is a global semiflow and it has a global attractor $\Cal A$. $\Cal A$ lies in $X^{1-\alpha}$ and is compact in the norm of $X^{1-\alpha}$. \endproclaim

\head\dfs{100405-1853}Preliminaries
\endhead
We assume the reader's familiarity with attractor theory on metric spaces as expounded
in e.g.~\cite{\rfa{Ha}} or, more recently,
in~\cite{\rfa{DCh}} and we just collect here a few relevant
concepts from that theory.
\definition{Definition}
Let $Y$ be a metric space.
Recall that a {\sl local semiflow $\pi$ on $Y$\/} is, by
definition, a continuous map from an open subset $D$ of
$\ro0,\infty..\times Y$ to $Y$ such that, for every $x\in
Y$ there is an $\omega_x=\omega_{\pi,x}\in\lo0,\infty..$
with the property that $(t,x)\in D$ if and only if $t\in
\ro0,\omega_x..$,  and such that (writing $x\pi
t:=\pi(t,x)$ for $(t,x)\in D$) $x\pi 0=x$ for $x\in Y$ and
whenever $(t,x)\in D$ and $(s,x\pi t)\in D$ then
$(t+s,x)\in D$ and $x\pi(t+s)=(x\pi t)\pi s$. Given an
interval $I$ in $\R$, a map $\sigma\co I\to Y$ is called a
{\sl solution (of $\pi$)\/} if whenever $t\in I$ and
$s\in\ro0,\infty..$ are such that $t+s\in I$, then
$\sigma(t)\pi s$ is defined and $\sigma(t)\pi
s=\sigma(t+s)$. If $I=\R$, then $\sigma$ is called a {\sl
full solution (of $\pi$)\/}. A subset $S$ of $Y$ is called
{\sl ($\pi$-)invariant\/} if for every $x\in S$ there is a
full solution $\sigma$ with $\sigma(\R)\subset S$ and
$\sigma(0)=x$. A point $x\in Y$ is called an {\sl equilibrium of $\pi$\/} if $x\pi t=x$ for all all $t\in \ro 0,\omega_x..$.

Given a local semiflow $\pi$ on $Y$ and a subset  $N$ of
$Y$, we say that {\sl $\pi$ does not explode in $N$\/} if
whenever $x\in Y$ and $x\pi \ro0,\omega_x..\subset N$, then
$\omega_x=\infty$. A {\sl global semiflow\/} is a local
semiflow with $\omega_x=\infty$ for all $x\in Y$.

Now let $\pi$ be a global semiflow on $Y$. A subset $A$ of
$Y$ is called a {\sl global attractor (rel. to $\pi$)\/} if
$A$ is  compact, invariant and if for every bounded set $B$
in $Y$ and every  open neighborhood $U$ of $A$ there is a
$t_{B,U}\in \ro0,\infty..$ such that $x\pi t\in U$ for all
$x\in B$ and all $t\in \ro t_{B,U},\infty..$. It easily
follows that a global attractor, if it exists, is uniquely
determined.

A subset $B$ of $Y$ is called {\sl ($\pi$-)ultimately
bounded} if there is a $t_B\in \ro0,\infty..$ such the set
$\{\,x\pi t\mid x\in B,\,t\in \ro t_B,\infty..\, \}$ is
bounded.

$\pi$ is called {\sl asymptotically compact} if whenever
$B\subset Y$ is  ultimately bounded, $(x_n)_n$ is a
sequence in $B$ and $(t_n)_n$ is a sequence in
$\ro0,\infty..$ with $t_n\to \infty$ as $n\to \infty$, then
the sequence $(x_n\pi t_n)_n$ has a convergent subsequence.
\enddefinition
The following result is well-known:
\proclaim{Proposition~\dft{250606-1446}} Let $\pi$ be a global semiflow on a metric space $Y$.
Suppose that
\roster
\item\dfro{270606-1419} $\pi$ is
asymptotically compact;
\item\dfro{270606-1420} every bounded subset of $Y$ is
ultimately bounded;
\item\dfro{270606-1421} the set of all equilibria of $\pi$ is bounded;
\item\dfro
{270606-1422} there is a continuous function $\Cal L\co Y\to \R$ which is bounded below, nonincreasing along solutions of $\pi$ and whenever $\Cal L(x\pi t)=\Cal L(x) $ for all $t\in\ro0,\infty..$ then $x$ is
an equilibrium of $\pi$.
\endroster
Under these assumptions, $\pi$ has a global attractor.
\endproclaim
\proof This is
just~\cite{\rfa{DCh}, Corollary~1.1.4 and
Proposition~1.1.3}.
\eproof

Given a Banach space $X$ and a sectorial operator $A\co D(A)\subset X\to X$ in $X$ with $\re\sigma(A)>0$ we know that $-A$ is the generator of an analytic semigroup $(e^{-At})_{t\in\ro0,\infty..}$ of linear operators on $X$. For $\alpha\in \oi0,\infty..$ we define, as usual, the operator $A^{-\alpha}\co X\to X$ as
$$A^{-\alpha}u=(1/\Gamma(\alpha))\int_0^\infty t^{\alpha-1}e^{-At}u\,\D t,\quad u\in X.\leqno\dff{230107-2127}$$
$A^{-\alpha}$ is injective and we define $X^\alpha=X_A^\alpha$ to be the range of $A^{-\alpha}$. We define $A^\alpha\co X^\alpha\to X$ to be the inverse of $A^{-\alpha}$. We also set $X^0=X$ and $A^{0}=\Id_X$. We call $A^{-\alpha}$, resp. $A^\alpha$ the {\it basic fractional power of $A$ of order $-\alpha$\/}, resp. $\alpha$ and we call $X^\alpha$ the {\it  fractional power space of $A$ of order $\alpha$\/}. $X^\alpha$ is a Banach space with respect to the norm
$$|u|_{X^\alpha}:=|A^\alpha u|_X,\quad u\in X^\alpha.$$ If $\beta>\alpha$ then $X^\beta$ is a dense subset of $X^\alpha$.

Moreover,
$$A^{-\beta}A^{-\gamma}x=A^{-\beta-\gamma}x,\quad \alpha,\beta\in \oi0,\infty..,\,x\in X.\leqno\dff{210107-2125}$$
Now let $X$ be a Hilbert space and $A\co D(A)\subset X\to X$ be self-adjoint in $X$ with $\re\sigma(A)>0$. Then $A$ is sectorial in $X$ and, for $\alpha\in\oi0,\infty..$,
$$A^{-\alpha}=\int_0^\infty t^{-\alpha}\,\D E(t),\leqno\dff{220107-1030}$$
where $(E(t))_{t\in\R}$ is the spectral measure defined by $A$.
In this case the set $X^\alpha$ is a Hilbert space with respect to the scalar product
$$\langle u,v\rangle_{X^\alpha}:=\langle A^\alpha u,A^\alpha v\rangle_X,\quad u,v\in X^\alpha. $$
For $\alpha\in\oi0,\infty..$ let
$X^{-\alpha}=X^{-\alpha}_A$ be the dual space of $X^\alpha$. We endow $X^{-\alpha}$ with the scalar product  $\langle\cdot,\cdot\rangle_{X^{-\alpha}}$ dual to the scalar product $\langle\cdot,\cdot\rangle_{X^\alpha}$, i.e.
$$\langle u',v'\rangle_{X^{-\alpha}}=\langle R_\alpha^{-1}u',R_\alpha^{-1}v'\rangle_{X^\alpha}, \quad u',v'\in X^{-\alpha},$$ where $R_\alpha\co X^\alpha\to X^{-\alpha}$ is the Fr\'echet-Riesz isomorphism $u\mapsto \langle\cdot,u\rangle_{X^\alpha}$.

 $X^{-\alpha}$ is called the {\it  fractional power space of $A$ of order $-\alpha$.\/}

By  $\phi$ denote
the duality map from $X$ to $X'$, i.e.
$$\phi(x):=\langle\cdot,x\rangle,\quad x\in
X.$$
Let $\alpha$, $\beta\in \R$ be arbitrary.  If $\beta\ge
\alpha\ge 0$, then let $\phi_{\beta,\alpha}\co X^\beta\to X^\alpha$ be
the inclusion map; if If $\beta\ge
\alpha> 0$, then define the map $\phi_{-\alpha,-\beta}\co X^{-\alpha}\to X^{-\beta}$ by
$\phi_{-\alpha,-\beta}(y')=y'_{|X^\beta}$ for $y'\in X^{-\alpha}$ i.e. for $y'\co X^\alpha\to \R$ linear and bounded; if $\beta>0$, define $\phi_{0,-\beta}\co X^0=X\to X^{-\beta}$ as follows: if $x\in X$, then $\phi_{0,-\beta}(x)$ is equal to the  map $y'\co X^\beta\to\R$ such that $y'(y)=\langle y,x\rangle$ for all $y\in X^\beta$. It follows that $y'=\phi_{0,-\beta}(x)\in X^{-\beta}$, so $\phi_{0,-\beta}$ is defined. Finally,
 if $\alpha>0$ and $\beta>0$, then let
$\phi_{\beta,-\alpha}:=\phi_{0,-\alpha}\circ\phi_{\beta,0}$.

We have the following basic result.
\proclaim{Proposition~\dft{280107-2215}}
For all $\alpha$, $\beta\in \R$ with $\beta\ge \alpha$ the map $\phi_{\beta,\alpha}\co X^\beta\to X^\alpha$ is defined, linear, bounded and injective. The set $\phi_{\beta,\alpha}[X^\beta]$ is dense in the Hilbert space $X^\alpha$. Moreover,
$$ \phi_{\alpha,\alpha}=\Id_{X^\alpha},\quad  \alpha\in\R$$
and
$$ \phi_{\gamma,\alpha}=\phi_{\beta,\alpha}\circ\phi_{\gamma,\beta}
, \quad\alpha,\beta,\gamma\in\R,\,
\gamma\ge\beta\ge\alpha.$$
For all $\alpha$, $\gamma\in\R$,  $\theta\in\ci0,1..$ and $x\in X^\gamma$ with $\alpha\le \gamma$ and $\beta=(1-\theta)\alpha+\theta \gamma$, the interpolation inequality
$$|\phi_{\gamma,\beta}x|_{X^\beta}\le |\phi_{\gamma,\alpha} x|_{X^\alpha}^{1-\theta}|x|_{X^\gamma}^\theta$$ holds.

For all $\alpha$, $\beta\in \ro0,\infty..$,
 $A_{(\alpha)}^{-\beta}:=A^{-\beta}|X^\alpha\co X^\alpha\to X^{\beta+\alpha}$  is a linear bijective isometry.

 For every $\alpha\in \oi0,\infty..$,
$\beta\in
\ro0,\infty..$  there is a unique continuous map
$A_{(-\alpha)}^{-\beta}\co X^{-\alpha}\to
X^{\beta-\alpha}$ with $A_{(-\alpha)}^{-\beta}\circ
\phi_{0,-\alpha}=\phi_{\beta,\beta-\alpha}\circ
A^{-\beta}$. $A_{(-\alpha)}^{-\beta}$ is a linear bijective isometry.

For $\alpha\in\R$ and $\beta\in\oi0,\infty..$  define the map
$A_{(\alpha)}^\beta\co X^\alpha\to X^{-\beta+\alpha}$ by
$$A_{(\alpha)}^\beta=(A_{(-\beta+\alpha)}^{-\beta})^{-1} $$
and set $A_{(\alpha)}:=A_{(\alpha)}^1$.
Then
for all $\gamma$, $\gamma'\in\R$ with $\gamma>\gamma'$ and
all $\beta\in\R$,
$$\phi_{-\beta+\gamma,-\beta+\gamma'}\circ
A_{(\gamma)}^{\beta}=A_{(\gamma')}^{\beta}\circ
\phi_{\gamma,\gamma'}$$
and
for all $\alpha$, $\beta$ and $\gamma\in\R$,
$$A_{(-\gamma+\alpha)}^\beta\circ A_{(\alpha)}^\gamma=A_{(\alpha)}^{\beta+\gamma}. $$
For $\alpha$, $\beta\in\R$ with $\beta\ge \alpha$ the map $\phi_{\beta,\alpha}$ is bijective from $ X^\beta$ to $\phi_{\beta,\alpha}[X^\beta]$. For $\alpha$, $\beta\in\oi0,\infty..$ define the map
$$\widetilde A_{(-\alpha)}^\beta:=A_{(\beta-\alpha)}^\beta\circ \phi_{\beta-\alpha,-\alpha}^{-1}\co \phi_{\beta-\alpha,-\alpha}[X^{\beta-\alpha}]\subset X^{-\alpha}\to X^{-\alpha}$$ and set $\widetilde A_{(-\alpha)}:=\widetilde A_{(-\alpha)}^1$. The map
 $\widetilde A_{(-\alpha)}^\beta$ is bijective and its inverse is
 $\widetilde
A_{(-\alpha)}^{-\beta}:=\phi_{\beta-\alpha,-\alpha}\circ A_{(-\alpha)}^{-\beta}$.

 For every $\alpha\in\oi0,\infty..$ the map
$B:=\widetilde A_{(-\alpha)}\co D(B)=\phi_{1-\alpha,-\alpha}[X^{1-\alpha}]\subset X^{-\alpha}\to X^{-\alpha}$ is self-adjoint in $X^{-\alpha}$ and $\re\sigma(B)>0$. For $\beta\in\oi0,\infty..$ let  $B^{-\beta}$ be the basic fractional power of $B$ of order $-\beta$ and $X^\beta_B$ be the corresponding fractional power space. Then
$$B^{-\beta}=\widetilde A_{(-\alpha)}^{-\beta}$$ and
$$X^\beta_B=\phi_{\beta-\alpha,-\alpha}[X^{\beta-\alpha}].$$
The map $\phi_{\beta-\alpha,-\alpha}$ is an isometry of the Hilbert space $X^{\beta-\alpha}$ onto $X^\beta_B$.

Finally, whenever $\alpha\in\ro0,(1/2)..$, $x\in X^{1-\alpha}$ and $v\in X^{1/2}\subset X^\alpha$, then $$(A_{(1-\alpha)}x).v=\langle  x,v\rangle_{X^{1/2}}.$$
Here, the dot `$.$' denotes function application between an element of $X^{-\alpha}$ and $X^\alpha$.
\endproclaim
\remark{Remark~\dft{220107-1421}} In view of  Proposition~\rf{280107-2215}, for $\alpha$, $\beta\in \R$ with $\beta\ge \alpha$ one often  regards $\phi_{\beta,\alpha}$ as an inclusion map and  $X^\beta$ as a (dense) subset of $X^\alpha$.

Sometimes (cf e.g.~\cite{\rfa{PR2}}) the notation
$$H_\alpha:=X^{\alpha/2}, \quad \alpha\in\R$$ is used. We then set
$$A_\alpha^\beta:=A_{(\alpha/2)}^\beta,\quad \alpha,\,\beta\in\R$$
and
$$A_\alpha:=A_\alpha^1,\quad \alpha\in\R.$$
Notice that $A_\alpha\co H_\alpha\to H_{\alpha-2}$ for all $\alpha\in\R$.
\endremark
\medskip\par
Proposition~\rf{280107-2215} is well-known (see e.g. the book of Amann \cite{\rfa{A}}) but it is not easy to find in the literature a proof that is both elementary and complete.  Therefore, in the Appendix, we provide  an elementary proof which presupposes only minimal knowledge of spectral measures.

\head\dfs{generalresults}Some results on semilinear parabolic equations\endhead

\proclaim{Proposition~\dft{150107-1446}}
Let $X$ be a Banach space and $A\co D(A)\subset X\to X$ be sectorial. Let $\overline k\in \ro0,\infty..$  be such that $\re\sigma(A+\overline kI)>0$ and $X^\beta$, $\beta\in \ro0,\infty..$, be the family of fractional power spaces generated by $A+\overline kI$. Let $\alpha\in \ro0,1..$, $g\co X^\alpha\to X$ be Lipschitzian on bounded subsets of $X^\alpha$ and $\pi$ be the local semiflow on $X^\alpha$ generated by the solutions of the differential equation
$$\dot u+Au=g(u).\leqno\dff{150107-1448}$$
Suppose $T\in \oi0,\infty..$, $u\co \ci0,T..\to X^\alpha$ and $u_k\co \ci0,T..\to X^\alpha$, $k\in \N$, are solutions of $\pi$ such that, for some $R\in\ro0,\infty..$, $|u(t)|_{X^\alpha}\le R$ and $|u_k(t)|_{X^\alpha}\le R$ for all $k\in \N$ and $t\in \ci0,T..$. If  $u_k(0)\to u(0)$ in $X$ as $k\to \infty$, then, for every $T_0\in \lo0,T..$, $u_k(t)\to u(t)$ in $X^\alpha$ as $k\to \infty$, uniformly for $t\in \ci T_0,T..$.
\endproclaim
\proof
There is a constant $L=L(R)\in\ro0,\infty..$ such that $$|g(v_1)-g(v_2)|_{X}\le L|v_1-v_2|_{X^\alpha}$$ for all $v_1$, $v_2\in X^\alpha$ with $|v_1|_{X^\alpha}\le R$ and $|v_2|_{X^\alpha}\le R$. Now, for every $k\in \N$ and $t\in \lo0,T..$,
$$u_k(t)-u(t)=e^{-At}(u_k(0)-u(0))+\int_0^t e^{-A(t-s)}(g(u_k(s))-g(u(s)))\,\D s,$$ so
$$\aligned|u_k(t)-u(t)|_{X^\alpha}&\le Ct^{-\alpha}|u_k(0)-u(0)|_X+C\int_0^t (t-s)^{-\alpha}|g(u_k(s))-g(u(s))|_X\,\D s\\
&\le Ct^{-\alpha}|u_k(0)-u(0)|_X+CL\int_0^t (t-s)^{-\alpha}|u_k(s)-u(s)|_{X^\alpha}\,\D s,\endaligned$$
for some constant $C\in\ro0,\infty..$, depending only on $\alpha$.
By Henry's inequality, cf.~\cite{\rfa{He}, Theorem 7.1.1} or~\cite{\rfa{DCh}, Lemma 1.2.9}, this implies that
$$|u_k(t)-u(t)|_{X^\alpha}\le C' t^{-\alpha}|u_k(0)-u(0)|_X,\quad k\in\N,\,t\in\lo0,T...\leqno\dff{150107-1910}$$
where $C'\in\ro0,\infty..$ is a constant which only depends on $(\alpha, C,L,T)$. Estimate~\rf{150107-1910} implies the assertion of the Proposition.
\eproof
\proclaim{Theorem~\dft{291206-1641}}
Let $X$ be a Hilbert space and $A\co D(A)\subset X\to X$ be selfadjoint and bounded from below. Let $\overline k\in \ro0,\infty..$  be such that $\re\sigma(A+\overline kI)>0$ and $X^\beta$, $\beta\in \R$, be the family of fractional power spaces generated by $A+\overline kI$. Let $\alpha\in \ro0,1..$, $g\co X^\alpha\to X$ be Lipschitzian on bounded subsets of $X^\alpha$ and $\pi$ be the local semiflow on $X^\alpha$ generated by the solutions of the differential equation
$$\dot u+Au=g(u).\leqno\dff{291206-1651}$$ If $K\subset X^\alpha$ is a $\pi$-invariant set which is compact in $X^\alpha$ then
$K\subset X^1=D(A)$ and $K$ is compact in $X^1$.

\endproclaim
\proof
By results in~\cite{\rfa{He}}, $K\subset X^1$. Let $(\overline u_n)_n$ be an arbitrary sequence in $K$. Then there is a sequence $(u_n)_n$ of solutions of $\pi$ lying in $K$ such that $u_n(0)=\overline u_n$ for every $n\in \N$. Let $\beta\in \oi0,1..$ be such that $\beta>\alpha$. By~\cite{\rfa{He}, Theorem 3.5.2, and its proof} there is a constant $C\in\oi0,\infty..$ such that for every $n\in\N$, $u_n$ is differentiable into $X^\beta$ and $$|v_n|_\beta\le C,\quad n\in\N\leqno\dff{291206-2002}$$ where $v_n=\partial(u_n;X^\beta)(0)$ for $n\in\N$. There is a strictly increasing sequence $(n_m)_m$ in $\N$ and a $\overline u\in K$ such that $\overline u_{n_m}\to \overline u$ in $X^\alpha$ as $m\to \infty$. Thus, using the notation of Proposition~\rf{280107-2215}, $\phi_{0,\alpha-1}(A+kI)\overline u_{n_m}\to \phi_{0,\alpha-1}(A+kI)\overline u$ in $X^{\alpha-1}$ and $g(\overline u_{n_m})\to g(\overline u)$ in $X^0$ as $m\to \infty$ so $-\phi_{0,\alpha-1}A\overline u_{n_m}+\phi_{0,\alpha-1}g(\overline u_{n_m})\to -\phi_{0,\alpha-1}A\overline u+\phi_{0,\alpha-1}g(\overline u)$ in $X^{\alpha-1}$ as $m\to\infty$. Now
$$\partial (u_n;X^\alpha)(t)=-Au_n(t)+g(u_n(t)),\quad n\in\N,\,t\in\R.$$ We thus conclude that
$\phi_{0,\alpha-1}v_{n_m}\to -\phi_{0,\alpha-1}A\overline u+\phi_{0,\alpha-1}g(\overline u)$ in $X^{\alpha-1}$ as $m\to \infty$. This together with~\rf{291206-2002} and the interpolation inequality  from Proposition~\rf{280107-2215} implies that $v_{n_m}\to -A\overline u+g(\overline u)$ in $X^{0}$ as $m\to \infty$. Thus $-A\overline u_{n_m}+g(\overline u_{n_m})\to -A\overline u+g(\overline u)$ in $X^{0}$ as $m\to\infty$ and as $g(\overline u_{n_m})\to g(\overline u)$ in $X^0$ as $m\to \infty$ it follows that $A\overline u_{n_m}\to A\overline u$ in $X^{0}$ as $m\to\infty$. It follows that $(A+kI)\overline u_{n_m}\to (A+kI)\overline u$ in $X^{0}$ so $\overline u_{n_m}\to \overline u$ in $X^{1}$. The theorem is proved.
\eproof

\head\dfs{linearestimates}Some linear estimates\endhead

\remark{Remark~\dft{270206-0701}}
 Under Hypothesis~\rf{221105-1756} item~$(1)$ let the
 operator $\Li\co H^1_0(\Omega)\to{\Cal D}'(\Omega)$ be
 defined by
 $$
  \Li u=\sum_{i,j=1}^N\partial_i(a_{ij}\partial_ju),\quad u\in H^1_0(\Omega).
 $$
 The definition of distributional derivatives implies that
 $$
  (\Li u-\beta u)(v)=-\int_\Omega[\sum_{i,j=1}^N a_{ij}\partial_i u\partial _j v+\beta uv]\,\D x,\quad u\in H^1_0(\Omega),\,v\in \Cal D(\Omega).
 \leqno\dff{260206}
 $$
 It follows by density that
 $$\mycenter{
  $\langle (\Li u-\beta u),v\rangle_{L^2}=-\displaystyle\int_\Omega[\sum_{i,j=1}^N a_{ij}\partial_i u\partial _j v+\beta uv]\,\D x$ for $u$, $v\in H^1_0(\Omega)$ with $\Li u-\beta u\in L^2(\Omega)$.
 }\leqno\dff{270206-0706} $$
\endremark \proclaim{Lemma~\dft{111105-1838}} Assume
Hypothesis~\rf{221105-1756}. If $\kappa\in\ro0,\lambda_1..$
is arbitrary and if $\overline\eps$ and $\rho$ are chosen
such that
  $\overline \eps\in\oi0,a_0..$, $\rho\in\oi0,1..$
and $c:=\min
\bigl(\rho(a_0-\overline\eps),(1-\rho)(\lambda_1-\kappa)-\rho(\overline
\eps+C_{\overline\eps}+\kappa )\bigr)>0$ then $$\aligned c(|\nabla
u|_{L^2}^2+|u|_{L^2}^2)&\le \int_\Omega[\sum_{i,j=1}^N a_{ij}\partial_i u\partial _j u+(\beta-\kappa) |u|^2]\,\D x\\&\le C(|\nabla
u|_{L^2}^2+|u|_{L^2}^2),\,\, u\in H^1_0(\Omega)\endaligned$$ where
$C:=\max(a_1+\overline\eps,\overline\eps+C_{\overline\eps})$.
\endproclaim \proof This is just a simple
computation.\eproof
\proclaim{Lemma~\dft{270206-0836}}
 Assume Hypothesis~\rf{221105-1756}. For $u$, $v\in H^1_0(\Om)$ define
 $$\langle u,v\rangle_1=\int_\Omega[\sum_{i,j=1}^N a_{ij}\partial_i u\partial _j v+\beta uv]\,\D x .\leqno\dff{080505-1517}$$
 $\langle \cdot,\cdot\rangle_1$ is a scalar product on $H^1_0(\Omega)$ and the norm defined by this
 scalar product is equivalent to the usual norm on
 $H^1_0(\Om)$.

\endproclaim
\proof
 This follows from  Lemma~\rf{111105-1838}.
\eproof
\proclaim{Lemma~\dft{220506-0819}}
Suppose $(Y,\langle\cdot,\cdot\rangle_{Y})$ and $(X,\langle\cdot,\cdot\rangle_{X})$ are (real or complex)
Hilbert spaces such that $Y\subset X$, $Y$ is dense in $(X,\langle\cdot,\cdot\rangle_{X})$ and the inclusion
$(Y,\langle\cdot,\cdot\rangle_{Y})\to(X,\langle\cdot,\cdot\rangle_{X})$
is continuous.
Then for every $u\in X$ there exists a unique $w_u\in Y$ such that
$$\langle v,w_u\rangle_Y=\langle v,u\rangle_X\text{ for all $v\in Y$.}$$
The map $B\co X\to X$, $u\mapsto w_u$ is linear, symmetric and positive. Let $B^{1/2}$ be a square root of $B$, i.e. $B^{1/2}\co X\to X$ linear, symmetric and $B^{1/2}\circ B^{1/2}=B$.
Then   $B$ and $B^{1/2}$ are injective and $R(B)$ is dense in $Y$. Set $X^{1/2}=X^{1/2}_B=R(B^{1/2})$ and $B^{-1/2}\co X^{1/2}\to X$ be the inverse of $B^{1/2}$. On $X^{1/2}$ the assignment $\langle u,v\rangle_{1/2}:=\langle B^{-1/2}u, B^{-1/2}v\rangle_X$ is a complete scalar product. We have $Y=X^{1/2}$ and $\langle \cdot,\cdot\rangle_{Y}=\langle \cdot,\cdot\rangle_{1/2}$.
\endproclaim
\proof The function $v\mapsto \langle v, u\rangle_X$ is linear and continuous on $Y$. Thus Fr\'echet-Riesz theorem implies the existence and uniqueness of $w$ and the linearity of $B$. Since, for $u$ and $v\in X$  $$\langle Bu,v\rangle_X=\langle Bu,Bv\rangle_Y=\langle u,Bv\rangle_X\leqno\dff{220506-0917}$$ it follows that $B$ is symmetric and positive. If $u\in X$ and $Bu=0$ then $0=\langle v,Bu\rangle_Y=\langle v,u\rangle_X$ for all $v\in Y$ and since $Y$ is dense in $X$ we see that $u=0$ so $B$ is injective. It follows that $B^{1/2}$ is injective as well. If $v\in Y$ and $\langle v,Bu\rangle_Y=0$ for all $u\in X$ then $\langle v,u\rangle_X=0$ for all $u\in X$ so $v=0$. It follows that $R(B)$ is dense in $Y$.
Clearly $\langle\cdot,\cdot\rangle_{1/2}$ is a complete scalar product on $X^{1/2}$. If $v\in X^{1/2}$ and $u\in X$ then
$$\langle v, Bu\rangle_{1/2}=\langle B^{-1/2}v,B^{1/2}u\rangle_X=
\langle v,u\rangle_X.$$
Thus if $\langle v,Bu\rangle_{1/2}=0$ for all $u\in X$, then $v=0$. This shows that $R(B)$ is dense in $X^{1/2}$.
We claim that
$$\langle u,v\rangle_Y=\langle u,v\rangle_{1/2}, \quad u,v\in R(B). \leqno\dff{220506-1035}$$
In fact, if $u$ and $v\in R(B)$ then $u=B\tilde u$ and $v=B\tilde v$ for some $u$ and $v\in X$. Thus
$$\langle u,v\rangle_Y=\langle u,B\tilde v\rangle_Y=\langle u,\tilde v\rangle_X$$
and
$$\langle u,v\rangle_{1/2}=\langle B^{1/2}\tilde u, B^{1/2}\tilde v\rangle_X=\langle B\tilde u, \tilde v\rangle_X=\langle u,\tilde v\rangle_X.$$
The claim is proved.

Since
$B^{1/2}$ is continuous from $X$ to $X$ with bound $|B^{1/2}|$
it follows that, for all $u\in X^{1/2}$,
$$|u|_X\le |B^{1/2}|\,|B^{-1/2}u|_X =|B^{1/2}|\,|u|_{1/2}$$
so the inclusion map $(X^{1/2},|\cdot|_{1/2})\to (X,|\cdot|_X)$ is continuous.

Now, if $u\in Y$, then there is a sequence $(u_n)_n$ in $R(B)$ converging to $u$ in $Y$. It follows that $(u_n)_n$ is a Cauchy sequence in $Y$, so, by~\rf{220506-1035}, it is a Cauchy sequence in $X^{1/2}$ and so it converges to a $v\in X^{1/2}$. By what we have proved so far, $(u_n)_n$ converges to $v$ and to $u$ in $X$. Thus $u=v$ so $u\in X^{1/2}$. It follows that $Y\subset X^{1/2}$. The same argument, with `$Y$' and `$X^{1/2}$' exchanged with each other, proves that $X^{1/2}\subset Y$. The last statement of the lemma follows from~\rf{220506-1035} by density.
\eproof
\proclaim{Proposition~\dft{020706-1432}}
Let $D(\A)$ be the set of all $u\in H^1_0(\Omega)$ such that $Lu-\beta u\in L^2(\Omega)$. For $u\in D(\A)$ define
$$\A u=-Lu +\beta u.$$
Then $\A\co D(\A)\to L^2(\Omega)$, selfadjoint in $X=L^2(\Omega)$ with $\re\sigma(\A)>0$. Moreover, if $X^\alpha$, $\alpha\ge 0$, is the family of fractional power spaces generated by $\A$, then $X^{1/2}=H^1_0(\Omega)$ and the scalar product on $X^{1/2}$ is identical to the scalar product $\langle\cdot,\cdot\rangle_1$ defined in Lemma~\rf{270206-0836}.
\endproclaim
\proof Let $\langle\cdot,\cdot\rangle$ denote the scalar product of $L^2(\Omega)$.
From~\rf{270206-0706} we conclude that
$$\langle -\A u,v\rangle=\langle v,-\A u\rangle, \quad u,v\in D(-\A).\leqno\dff{020706-1825}$$
Lemma~\rf{111105-1838} implies that
$$\langle -\A u,u\rangle\le0,\quad u\in D(-\A).$$
We have thus proved that $-A$ is symmetric and dissipative.
We will now prove that $-\A$ is $m$-dissipative. To this end, we must prove that for every $\lambda\in \oi0,\infty..$ and every $g\in L^2(\Omega)$ there is a $u\in D(-\A)$ such that
$$u+\lambda \A u=g.$$
Define the bilinear form $b\co H^1_0(\Omega)\times H^1_0(\Omega)\to \R$ by
$$b(u,v)=\int_\Omega uv\,\D x+\lambda \int_\Omega[\sum_{i,j=1}^N a_{ij}\partial_i u\partial _j v+\beta uv]\,\D x,\quad u,v\in H^1_0(\Omega).$$
It follows from Hypothesis~\rf{221105-1756} and Lemma~\rf{111105-1838} that there are constants $C$ and $c\in \oi 0,\infty..$ such that, for $u$, $v\in H^1_0(\Omega)$
$$|b(u,v)|\le C|u|_{H^1_0}|v|_{H^1_0} $$
and
$$b(u,u)\ge c|u|^2_{H^1_0}.$$
Thus Lax-Milgram theorem shows that for every $g\in L^2(\Omega)$ there is a $u\in H^1_0(\Omega)$ such that $$b(u,v)=\langle g,v\rangle, \quad v\in H^1_0(\Omega).$$
In particular,
$u+\lambda(-Lu +\beta u)=g$ in the distributional sense.
It follows that $-Lu +\beta u\in L^2(\Omega)$ so $u\in D(\A)=D(-\A)$ and $u+\lambda \A u=g$. Therefore, indeed, $-\A$ is $m$-dissipative. Now an application of the results of~\cite{\rfa{CH}, Section 2.4} shows that $-\A$ is selfadjoint. Thus $\A$ is selfadjoint and $\re\sigma(\A)>0$ by Lemma~\rf{111105-1838}.
To prove the last statement of the Proposition,
set $(X,\langle\cdot,\cdot\rangle_X)=(L^2(\Omega),\langle \cdot,\cdot\rangle)$ and $(Y,\langle\cdot,\cdot\rangle_Y)=(H^1_0(\Omega),\langle \cdot,\cdot\rangle_1$, where the scalar product $\langle\cdot,\cdot\rangle_1$ is defined in Lemma~\rf{270206-0836}.
 Then $Y$ is dense in $X$ and the inclusion $Y\to X$ is continuous. Let $\bold B\co X\to X$ be the inverse of ${\bold A}$.
Then for all $u\in X$, $\bold Bu\in Y$ and \rf{270206-0706} implies that
 for all $v\in Y$
$$\langle v,u \rangle_X=\langle v,\bold Bu\rangle_Y.$$ Thus $\bold B=B$ where $B$ is as in Lemma~\rf{220506-0819}. Now Lemma~\rf{220506-0819} and Lemma~\rf{270206-0836} imply the proposition.
\eproof

\head\dfs{nonlinearestimates}Some nonlinear estimates\endhead

In this section we assume that $f\co\Omega\times\R\to\R$ is a function
satisfying Hypothesis \rf{231105-0903}.
\proclaim{Lemma~\dft{100107-0749}} Let $X$ be a Banach space and
$A\co D(A)\subset X\to X$ be a sectorial operator with
$\re\sigma(A)>0$ generating the basic family $X^\alpha$,
$\alpha\in\ro0,\infty..$ of fractional power spaces. Suppose that
$X^{1/2}$ is continuously included in $L^6(\Omega)$ and $X=X^0$ is
continuously included in $L^2(\Omega)$. Then for every
$p\in\ro2,6..$ there is a $\overline\beta\in \ro0,1/2..$ such that
for all $\alpha\in\oi\overline\beta,1..$ the space $X^\alpha$ is
continuously included in $L^p(\Omega)$.
\endproclaim
\proof If $u\in X^{1/2}$ then $u\in L^2(\Omega)\cap L^6(\Omega)$
so, by interpolation of Lebesgue spaces, $u\in L^p(\Omega)$ and
$$|u|_{L^p}\le |u|_{L^2}^{\beta}|u|_{L^6}^{1-\beta}$$ where
$\beta=(6-p)/(2p)$. Let $B$ be the inclusion map from $X^{1/2}$ to
$Y=L^p(\Omega)$. Since, by interpolation of fractional power
spaces, $$|u|_{X^{1/2}}\le C|u|_{X^0}^{1/2}|u|_{X^1}^{1/2}, \quad
u\in X^1$$ for some constant $C\in\ro0,\infty..$ we see that, for
$u\in D(A)=X^1\subset X^{1/2}=D(B)$, $$\aligned|Bu|_Y&\le
C'|u|_{X^0}^\beta|u|_{X^{1/2}}^{1-\beta}\le
C'|u|_{X^0}^\beta(C|u|_{X^0}^{1/2}|u|_{X^1}^{1/2})^{1-\beta}\\&=
C'C|u|_{X^1}^{\overline\beta}|u|_{X^0}^{1-\overline\beta}=
C'C|Au|_{X^0}^{\overline\beta}|u|_{X^0}^{1-\overline\beta}\endaligned$$
for some constant $C'\in\ro0,\infty..$, where
$\overline\beta=(1/2)(1-\beta)$. By~\cite{\rfa{He}, p. 28,
Exercise~11} we now obtain that for every $\alpha \in
\oi\overline\beta,1..$ the map $B\circ A^{-\alpha}$ is defined and
continuous from $X$ to $Y$. Thus $X^\alpha=R(A^{-\alpha})$ is
continuously included in $L^p(\Omega)$, as claimed. \eproof

\proclaim{Lemma~\dft{050207-1113}}
  Let $X$ be a Banach space and $A\co D(A)\subset X\to X$ be a sectorial operator with $\re\sigma(A)>0$,  generating the family
  $X^\alpha$, $\alpha\in\ro0,\infty..$ of fractional power spaces. Suppose that $X^{1/2}$ is continuously included in
  $H^1_0(\Omega)$ and $X=X^0$ is continuously included in $L^2(\Omega)$. Let $w\co\Omega\to\R$ be such that
  $w\in {\Cal E}_\gamma$ for some $\gamma\in]0,1[$.
 Then for all $\alpha\in\lo(\gamma/2),1..$, the mapping $u\mapsto |w|^{1/2}u$ defines a bounded linear fuction from $X^\alpha$ to $L^2(\Omega)$.
\endproclaim

\proof
It is easy to check that $w\in{\Cal E}_\gamma$ if and only if there exists a constant $C'$ such that
$$
||w|^{1/2}u|_{L^2}\leq C'|u|_{H^1}^\gamma|u|_{L^2}^{1-\gamma},\quad u\in H^1_0(\Omega).
$$
Let $B$ be the map from $X^{1/2}$ to $Y=L^2(\Omega)$ defined by the assignement $u\mapsto |w|^{1/2}u$. Since, by interpolation of fractional power spaces,
$$|u|_{X^{1/2}}\le C|u|_{X^0}^{1/2}|u|_{X^1}^{1/2}, \quad u\in X^1$$ for some constant $C\in\ro0,\infty..$ we see that, for $u\in D(A)=X^1\subset X^{1/2}=D(B)$,
$$\aligned|Bu|_Y&\le C'|u|_{X^{1/2}}^\gamma|u|_{X^{0}}^{1-\gamma}\le C'|u|_{X^{0}}^{1-\gamma}(C|u|_{X^{1}}^{1/2}|u|_{X^0}^{1/2})^{\gamma}\\&=
C'C|u|_{X^1}^{\gamma/2}|u|_{X^0}^{1-{\gamma/2}}= C'C|Au|_{X^0}^{\gamma/2}|u|_{X^0}^{1-{\gamma/2}}\endaligned$$
for some constant $C'\in\ro0,\infty..$. By~\cite{\rfa{He}, p. 28, Exercise~11} we now obtain that for every $\alpha \in \oi\gamma/2,1..$ the map $B\circ A^{-\alpha}$ is defined and continuous from $X$ to $Y$. Thus $B$ is a bounded linear map form  $X^\alpha=R(A^{-\alpha})$ to $Y=L^2(\Omega)$, as claimed.
\eproof

%\proclaim{Corollary~\dft{050207-1211}}
%Let $X$ and $A$ be as in Lemma~\rf{050207-1113} and assume that $X=L^2(\Omega)$.
%Let $w\co\Omega\to\R$ be such that $w\in {\Cal E}_\gamma$ for some $\gamma\in]0,1[$.
%Then for all $\alpha$, with $\gamma/2<\alpha\leq1$, the mapping $u\mapsto wu$ defines
%a bounded linear fuction from $X^\{\alpha}$ to $X^{-\alpha}$.
%\endproclaim

\proclaim{Proposition~\dft{100405-0808}}
Let $X=L^2(\Omega)$, $A\co D(A)\subset X\to X$ a sectorial operator with $\re\sigma(A)>0$,  generating the family $X^\alpha$, $\alpha\in\ro0,\infty..$ of fractional power spaces. Suppose that $X^{1/2}$ is continuously included in $H^1_0(\Omega)$. Let $q=(6/(\overline \rho+1))$ and $p=(q/(q-1))$.
 If $\overline \rho>2$ or $a^2\not\in{\Cal E}_1$, then choose $\alpha\in\oi0,(1/2)..$ such that
$\alpha>\max\{\gamma/2,(1-(6-p)/2p))/2\}$, so
$X^\alpha$ is continuously included in $L^p(\Omega)$ (by Lemma~\rf{100107-0749})
and the mapping $u\mapsto |a|^{1/2}u$ is bounded from $X^\alpha$ to $X$ (by Lemma~\rf{050207-1113}).
 If $\overline \rho= 2$ and $a^2\in{\Cal E}_1$ then let $\alpha=0$.

 Let $F\co \Om\times \R\to \R$ be defined by $$F(x,u)=\int_0^u
f(x,s)\,\D s,
 $$ whenever $s\mapsto f(x,s)$ is continuous and $F(x,u)=0$
 otherwise.
 Then $\widehat F$ maps $X^{1/2}$ into $L^1(\Omega)$ and the operator $\widehat F\co X^{1/2}\to L^1(\Omega)$ is Fr\'echet-differentiable with $D\widehat F(u).h=\widehat f(u)\cdot h$ for $u$, $h\in X^{1/2}$.

If $\overline\rho>2$,
  then for every $u\in X^{1/2}$ the function $\bold f(u)\co X^\alpha\to \R$, $$v\mapsto \int_\Omega  f(x,u(x))v(x)\,\D x,$$ is defined, linear and bounded, hence $\bold f(u)\in X^{-\alpha}$.

If $\overline\rho=2$, then  for every $u\in X^{1/2}$, let $\bold f(u):=\widehat f(u)\in X=X^{-\alpha}$.

For all $\overline \rho\in\ro2,4..$, the operator $\bold f\co X^{1/2}\to X^{-\alpha}$ is Lipschitzian on bounded subsets of $X^{1/2}$.
\endproclaim
\proof If $u\in X^{1/2}$ then $u\in L^2(\Omega)\cap L^6(\Omega)$ so $u\in L^r(\Omega)$ for every $r\in\ci2,6..$. In particular, $\widehat f(u)\co \Omega\to \R$ and $\widehat F(u)\co \Omega\to \R$ are measurable.

Now, for $u\in X^{1/2}$ we have that
$$\aligned\int_\Omega &|F(x,u(x))|\,\D x\\&\le \int_\Omega(|f(x,0)u(x)|+\overline C||a(x)|^{1/2}u(x)|^2/2+\overline C|u(x)|^{\overline \rho+2}/(\overline \rho+2))\,\D x<\infty\endaligned$$
as $\overline \rho+2\in\ci2,6..$. Hence $\widehat F(u)\in L^1(\Omega)$. Moreover, for $u$, $h\in X^{1/2}$ we have
$$\aligned &\int_\Omega|f(x,u(x))h(x)|\, \D x\\&\le\int_\Omega (|f(x,0)h(x)|+\overline C|a(x)u(x)h(x)|+\overline C|u(x)|^{\overline \rho+1}|h(x)|\,\D x)\\&\le
|\widehat f(0)|_{L^2}|h|_{L^2}+\overline C||a|^{1/2}u|_{L^2}||a|^{1/2}h|_{L^2}+\overline C|u|_{L^r}^{\overline\rho +1}|h|_{L^6},\endaligned$$
where $r=(6/5)(\overline \rho+1)$. It follows that for every $u\in X^{1/2}$ the map $h\mapsto \widehat f(u)\cdot h$ is linear and bounded from $X^{1/2}$ to $L^1(\Omega)$.
Now, for $u$, $h\in X^{1/2}$,
$$\aligned |\widehat F(u+h)&-\widehat F(u)-\widehat f(u)\cdot h|_{L^1}\\&=\int_\Omega|F(x,u(x)+h(x))-F(x,u(x))-f(x,u(x))h(x)|\,\D x \\&\le
\overline C
\int_\Omega(|a(x)h(x)h(x)|+\max(1,2^{\overline \rho-1})(|u(x)|^{\overline \rho}+|h(x)|^{\overline \rho})|h(x)||h(x)|)\,\D x\\&
\le \overline C||a|^{1/2}h|_{L^2}^2+\overline C \max(1,2^{\overline \rho-1})(|u|_{L^r}^{\overline \rho} +|h|_{L^r}^{\overline\rho})|h|_{L^6}^2,\endaligned$$
where $r=(6/4)\overline\rho$.
This shows that the operator $\widehat F\co X^{1/2}\to L^1(\Omega)$ is Fr\'echet-differentiable with $D\widehat F(u).h=\widehat f(u)\cdot h$ for $u$, $h\in X^{1/2}$.

Now suppose $\overline\rho>2$ or $a^2\not\in{\Cal E}_1$.
The fact that $X^\alpha$ is continuously embedded in $L^2(\Omega)$ and in $L^p(\Omega)$ (with a common embedding constant $C\in\ro0,\infty..$) implies that, for all $v\in X^\alpha$,
$$\aligned&\int_\Omega|f(x,u(x))v(x)|\,\D x\\&\le
\int_\Omega(|f(x,0)v(x)|+\overline C|a(x)u(x)v(x)|+\overline C|u(x)|^{\overline \rho+1}|v(x)|)\,\D x\\&\le
|\widehat f(0)|_{L^2}|v|_{L^2}+\overline C ||a|^{1/2}u|_{L^2}||a|^{1/2}v|_{L^2}+|u|_{L^6}^{\overline \rho+1}|v|_{L^p}\\&\le
C(|\widehat f(0)|_{L^2}+\overline C ||a|^{1/2}u|_{L^2}+|u|_{L^6}^{\overline \rho+1})|v|_{X^{\alpha}}.\endaligned$$
Thus, indeed, the function $\bold f(u)\co X^\alpha\to \R$, $v\mapsto \int_\Omega  f(x,u(x))v(x)\,\D x$, is defined, linear and bounded, hence $\bold f(u)\in X^{-\alpha}$. Similarly, we obtain for $u$, $h\in X^{1/2}$ and $v\in X^\alpha$,
$$\aligned&\int_\Omega|(f(x,u(x)+h(x))-f(x,u(x)))v(x)|\,\D x\\&\le \overline C
\int_\Omega(|a(x)h(x)v(x)|+\max(1,2^{\overline \rho-1})(|u(x)|^{\overline \rho}+|h(x)|^{\overline \rho})|h(x)||v(x)|)\,\D x\\&\le
\overline C||a|^{1/2}h|_{L^2}||a|^{1/2}v|_{L^2}+\overline C (|u|_{L^6}^{\overline \rho}+|h|_{L^6}^{\overline \rho})|h|_{L^6}|v|_{L^p}\\&\le
C(\overline C||a|^{1/2}h|_{L^2}+\overline C (|u|_{L^6}^{\overline \rho}+|h|_{L^6}^{\overline \rho})|h|_{L^6})|v|_{X^{\alpha}}.\endaligned$$
This shows that the operator $\bold f\co X^{1/2}\to X^{-\alpha}$ is defined and Lipschitzian on boun\-ded subsets of $X^{1/2}$.

Now suppose $\overline\rho=2$ and $a^2\in{\Cal E}_1$. Then similar arguments show that
$$|\widehat f(u)|_{L^2}\le |\widehat
f(0)|_{L^2}+\overline C(|a u|_{L^2}+ |u|^{\overline \rho
+1}_{L^{2(\overline \rho+1)}}), $$
$$\aligned|\widehat f(u+h)&-\widehat f(u)|_{L^2}\\&\le \overline
C|a h|_{L^2}+ \overline C\max(1,2^{\overline
\rho-1})(|u|^{\overline \rho}_
{L^{2(\overline\rho+1)}}+|h|^{\overline
\rho}_{L^{2(\overline\rho+1)}})
|h|_{L^{2(\overline\rho+1)}},
\endaligned
$$and
$$\aligned&|\widehat
F(u+h)-\widehat F(u)-\widehat f(u)h|_{L^1}\\&\le\bigl( \overline
C|a h|_{L^2}+ \overline C\max(1,2^{\overline
\rho-1})(|u|^{\overline \rho}_{L^{2(\overline\rho+1)}}+
|h|^{\overline
\rho}_{L^{2(\overline\rho+1)}})|h|_{L^{2(\overline\rho+1)}}\bigr)
|h|_{L^2}. \endaligned$$
Again this shows that the operator $\bold f\co X^{1/2}\to X^{-\alpha}=X$ is defined and Lipschitzian on bounded subsets of $X^{1/2}$.
\eproof

\head\dfs{100405-1955} Tail estimates and the existence of
attractors \endhead
Let $\A$ be the operator defined in Proposition~\rf{020706-1432} and $X^\alpha$, $\alpha\in\R$, $\A_{(\alpha)}$, $\alpha\in\R $  and $\widetilde\A_{(-\alpha)}$, $\alpha\in\oi0,\infty..$, be the spaces and the operators defined in Proposition~\rf{280107-2215} with respect to $A=\A$.

If $\overline \rho= 2$, then, by what we have proved so far, the parabolic equation
$$\dot u=-\A u+\bold f(u)$$ defines a local semiflow $\pi$ on $X^{1/2}$.

If $\overline\rho\in\oi2,4..$, then choose $\alpha\in\oi0,(1/2)..$ as in Proposition~\rf{100405-0808}. Then the parabolic equation
$$\dot{\widetilde u}=-\widetilde\A_{(-\alpha)}\widetilde u+\bold f(\phi_{(1/2),-\alpha}^{-1}\widetilde u)$$
defines a local semiflow $\widetilde \pi$ on $\phi_{(1/2),-\alpha}[X^{1/2}]$.

Let $\pi$ be the local semiflow on $X^{1/2}$ which is conjugate to $\widetilde \pi$ via the conjugation $\phi_{(1/2),-\alpha}\co X^{1/2}\to \phi_{(1/2),-\alpha}[X^{1/2}]$.

While the local semiflow $\widetilde \pi$ depends on $\alpha$, it is not difficult to prove that
\proclaim{Proposition~\dft{310107-1126}} For $\overline\rho\in\oi2,4..$, the local semiflow $\pi$ is independent of the choice of $\alpha\in\oi0,(1/2)..$ such that $X^\alpha$ is continuously imbedded in $L^p(\Omega)$, where $p=(q/(q-1))$ and $q=(6/(\overline\rho+1))$.

\endproclaim
From the definition of $\pi$ we thus obtain, choosing $\alpha=0$ if $\overline\rho= 2$, that
\proclaim{Proposition~\dft{310107-1256}}
Whenever $T\in\lo0,\infty..$ and $u\co \ro0,T..\to X^{1/2}$ is a solution of $\pi$, then $u$ is continuous on $\ro0,T..$, differentiable (into $X^{1/2}$) on $\oi0,T..$ and, for $t\in\oi0,T..$
$$\phi_{0,-\alpha}(\dot u(t))=-\A_{1-\alpha}u(t)+\bold f(u(t)),$$
where $\dot u(t):=\partial(u;X^{1/2})(t)=\partial(u;X^0)(t)$.
\endproclaim
\proclaim{Proposition~\dft{270606-1118}}
Define the function $\Cal L\co H^1_0(\Omega)\to \R$ by
$$\Cal L(u)=(1/2)\langle u,u\rangle_{X^{1/2}}-\int_{\Omega}F(x,u(x))\,\D x,\, u\in H^1_0(\Omega).$$
Let $T\in\lo
0,\infty..$  and $u\co \ro0,T..\to H^1_0(\Omega)=X^{1/2}$ be a
solution of  $\pi$.
 Then the function
$\Cal L\circ u$ is continuous on $\ro0,T..$,  differentiable on $\oi0,T..$ and, for
$t\in \oi0,T..$,
$$(\Cal L\circ u)'(t)=-|\dot u(t)|_{L^2}.\leqno\dff{310107-1331}$$
\endproclaim
\proof
By Proposition~\rf{100405-0808}  the function $\Cal L$ is Fr\'echet differentiable and
$$D\Cal L(w).v=\langle w,v\rangle_{X^{1/2}}-\int_{\Omega}\widehat f(w(x))v(x)\,\D x,\quad w,\,v\in X^{1/2}.$$
Thus, by Proposition~\rf{310107-1256}, $\Cal L\circ u$ is continuous on $\ro0,T..$,  differentiable on $\oi0,T..$ and, for
$t\in \oi0,T..$,
$$(\Cal L\circ u)'(t)=\langle u(t),\dot u(t)\rangle_{X^{1/2}}-\int_{\Omega}\widehat f(u(t)(x))\dot u(t)(x)\,\D x.$$
If $\overline\rho>2$, then the last statement of Proposition~\rf{280107-2215} implies that
$$\langle u(t),\dot u(t)\rangle_{X^{1/2}}=(\A_{1-\alpha}u(t)).\dot u(t)$$ and so, by Proposition~\rf{100405-0808} and~\rf{310107-1256},
$$\aligned(\Cal L\circ u)'(t)&=(\A_{1-\alpha}u(t)).\dot u(t)-\bold f(u(t)).\dot u(t)=(\A_{1-\alpha}u(t)-\bold f(u(t))).\dot u(t)\\&=(-\phi_{0,-\alpha}\dot u(t)).\dot u(t)=-\langle \dot u(t),\dot u(t)\rangle_{X}\endaligned$$
This proves formula~\rf{310107-1331} for $\overline\rho\in\oi2,4..$. A similar but simpler argument proves~\rf{310107-1331} for $\overline\rho=2$.
\eproof
\proclaim{Corollary~\dft{270606-1148}}
$\pi$ is a global semiflow on $Y=H^1_0(\Omega)$ and it satisfies properties~\rf{270606-1420}, \rf{270606-1421} and~\rf{270606-1422} of Proposition~\rf{250606-1446}.
\endproclaim
\proof
Proposition~\rf{270606-1118} together with Proposition~\rf{100405-0808} implies that $\Cal L$ is continuous and nonincreasing along solutions of $\pi$.
Thus, for all $u_0\in H^1_0(\Omega)$, writing $u(t)=u_0\pi t$ for $t\in\ro0,\omega_{u_0}..$, we obtain
$$\aligned(1/2)\langle u(t),u(t)\rangle_{X^{1/2}}&\le (1/2)\langle u_0,u_0\rangle_{X^{1/2}}-\int_\Omega F(x,u_0(x))\,\D x\\&+\int_\Omega F(x,u(t)(x))\,\D x\le
(1/2)\langle u_0,u_0\rangle_{X^{1/2}}-\int_\Omega F(x,u_0(x))\,\D x\\&+\int_\Omega c(x)\,\D x,\quad t\in\ro0,\omega_{u_0}...\endaligned\leqno\dff{270606-1437}$$
It follows that every solution of $\pi$ is bounded in $Y$ and so, as $\pi$ does not explode in bounded subsets of $Y$, $\pi$ is a global semiflow. Proposition~\rf{100405-0808} also implies that every bounded set $B$ is $\pi$-ultimately bounded, with $t_B=0$.
Moreover, $\Cal L (u)\ge -\int_\Omega c(x)\,\D x$ for every $u\in Y$ so $\Cal L$ is bounded from below. If $u\co \ro0,\infty..\to Y$ is a solution of $\pi$ along which $\Cal L$ is constant, then, by Proposition~\rf{270606-1118}, $\dot u(t)=0$ for all $t\in\oi0,\infty..$ and this clearly implies that $u$ is a constant solution so $u(0)$ is an equilibrium of $\pi$. Let $u_0$ be an equilibrium of $\pi$. Suppose first that $\overline\rho>2$. Then ${\bold A}_{(1-\alpha)}u_0$ is defined and ${\bold A}_{(1-\alpha)}u_0=\bold f(u_0)$. It follows that
$$\langle u_0,u_0\rangle_{X^{1/2}}= ({\bold A}_{(1-\alpha)}u_0).u_0=(\bold f(u_0)).u_0=\int_{\Omega}f(x,u_0(x))u_0(x)\,\D x\le\int_\Omega c(x)\,\D x$$ so the set of all equilibria of $\pi$ is bounded in $Y$. A similar but simpler argument shows that, if $\overline\rho= 2$, then again the set of all equilibria of $\pi$ is bounded in $Y$. The proposition is proved.
\eproof
We can now state our basic result on tail estimates.
\proclaim{Theorem~\dft{270606-1543}}
 $\overline\vartheta\co
\R\to \ci 0,1..$ be a function of class $C^1$ with $\overline\vartheta (s)=0$ for $s\in
\lo-\infty,1..$ and $\overline\vartheta (s)=1$ for $s\in
\ro2,\infty..$. Let $\vartheta:=\overline\vartheta^2$.
For $k\in\N$ let the functions
$\overline\vartheta_k\co\R^N\to\R$ and $\vartheta_k\co
\R^N\to\R$ be defined by
$$\overline\vartheta_k(x)=\overline\vartheta(|x|^2/k^2)\text{
and }\vartheta_k(x)=\vartheta(|x|^2/k^2),\quad x\in \R^N.$$
Set $C_\vartheta=2\sqrt
2\sup_{y\in\R}|\vartheta'(y)|$ and
$C_{\overline\vartheta}=2\sqrt
2\sup_{y\in\R}|\overline\vartheta'(y)|$.
Let $\kappa\in\oi0,\lambda_1..$ be arbitrary, where $\lambda_1$ is defined in Hypothesis~\rf{221105-1756}. For every $k\in\N$ let
$$b_k=\max(a_1C_{\overline \vartheta}^2k^{-2},2a_1C_{\overline\vartheta}k^{-1},a_1C_{\vartheta}k^{-1})$$ and $c_k=\int_\Omega \vartheta_k(x)c(x)\,\D x$.

Then whenever $R\in\ro0,\infty..$, $\tau\in\oi0,\infty..$ and $u\co\ro0,\infty..\to H^1_0(\Omega)$ is a solution of $\pi$ such that $|u(t)|_{H^1_0}\le R$ for all $t\in \ci0,\tau..$, then, for every $t\in\ci0,\tau..$,
$$\int_\Omega \vartheta_k(x)|u(x)|^2\,\D x\le R^2e^{-2\kappa t}+(b_kR^2+c_k)/\kappa.\leqno\dff{010207-0925}$$
\endproclaim

\proof
Since $\nabla
\vartheta_k(x)=(2/k^2)\vartheta'(|x|^2/k^2)x$ and $\nabla
\overline\vartheta_k(x)=(2/k^2)\overline\vartheta'(|x|^2/k^2)x$
we see that $$\text{$\sup_{x\in\Omega}|\nabla \vartheta_k
(x)|\le C_\vartheta /k$ and $\sup_{x\in\Omega}|\nabla
\overline\vartheta_k (x)|\le C_{\overline\vartheta}
/k$.}\leqno\dff{171105-0943}$$

Assume first that $\overline\rho>2$.

We claim that, for all $v\in H^1_0(\Omega)=X^{1/2}$,
$$\aligned(-\A_{(1-\alpha)}v)&.(\vartheta_k v)+\kappa\int_{\Omega}\vartheta_k(x)|v(x)|^2\,\D x\\&\le
\int_\Omega \bigl( (a_1C_{\overline \vartheta}^2k^{-2})|v|^2+(2a_1C_{\overline\vartheta}k^{-1})|v||\nabla v|+(a_1C_{\vartheta}k^{-1})|v||\nabla v|\bigr)\D x.\endaligned\leqno\dff{310107-1838}$$
To prove~\rf{310107-1838}, we may suppose that $v\in X^1$ since the general case follows by the density of $X^1$ in $H^1_0(\Omega)$.
Now, if $v\in X^1$, then
$$(-\A_{(1-\alpha)}v).(\vartheta_k v)=(-\phi_{0,-\alpha}\A v). (\vartheta_k v)=\langle -\A v,\vartheta_k v\rangle_X$$ so, using~\rf{171105-0943}, we obtain
$$\aligned &(-\A_{(1-\alpha)}v).(\vartheta_k v)+\kappa\int_{\Omega}\vartheta_k(x)|v(x)|^2\,\D x=\int_\Omega\vartheta_k v(x)(-{\bold A}v)(x)\,\D x\\&+\kappa\int_\Omega \vartheta_k|v(x)|^2\,\D x=
\int_\Omega\left(-(A\nabla v)(x)\cdot\nabla(\vartheta_k v)(x)
-\vartheta_k\beta(x)|v(x)|^2\right)\,\D x\\&+\kappa\int_\Omega \vartheta_k|v|^2\,\D x=
\int_\Omega\left(-(A\nabla(\overline \vartheta_kv))\cdot\nabla(\overline\vartheta_kv)
-(\beta(x)-\kappa)|\overline \vartheta_kv|^2\right)\,\D x\\&+
\int_\Omega\left(|v|^2(A\nabla\overline\vartheta_k) \cdot\nabla\overline\vartheta_k
+2v\overline\vartheta_k(A\nabla\overline\vartheta_k)\cdot\nabla v- v(A\nabla\vartheta_k)\cdot\nabla v\right)\,\D x\\
&\le\int_\Omega \left( (a_1C_{\overline \vartheta}^2k^{-2})|v|^2+(2a_1C_{\overline\vartheta}k^{-1})|v||\nabla v|+(a_1C_{\vartheta}k^{-1})|v||\nabla v|\right)\D x.\endaligned$$
This proves~\rf{310107-1838}.

For $k\in\N$ define the function $V_k\co H^1_0(\Omega)\to \R$
by
$$V_k(u)=(1/2)\int_\Omega \vartheta_k(x)|u(x)|^2\,\D x, \quad u\in H^1_0(\Omega).$$
Then $V_k$ is Fr\'echet differentiable and
$$DV_k(u)v=\int_\Omega \vartheta_k(x)u(x)v(x)\,\D x, \quad u,\,v\in H^1_0(\Omega).$$ Let $u\co \ro0,\infty..\to H^1_0(\Omega)$ be a solution of $\pi$. By Proposition~\rf{310107-1256},  $V_k\circ u$ is differentiable on $\oi0,\infty..$ and, for $t\in\oi0,\infty..$,
$$(V_k\circ u)'(t)=\int_{\Omega}\vartheta_k(x)u(t)(x)\dot u(t)(x)\,\D x=\langle \dot u(t),\vartheta_k u(t)\rangle_X.$$
By results in~\cite{\rfa{PR1}}, $\vartheta_ku(t)\in H^1_0(\Omega)=X^{1/2}\subset X^\alpha$. It follows that, for $t\in\oi0,\infty..$,
$$\aligned\langle \dot u(t),\vartheta_k u(t)\rangle_X&=(\phi_{0,-\alpha}\dot u(t)).(\vartheta_k u(t))=(-\A_{(1-\alpha)}u(t)+\bold f (u(t))).(\vartheta_k u(t))\endaligned$$
so, using~\rf{310107-1838}, we obtain
$$\aligned &(V_k\circ u)'(t)+2\kappa (V_k\circ u)(t)=(-\A_{(1-\alpha)}u(t)).(\vartheta_k u(t))\\&+\int_\Omega \widehat f(u(t))(x)(\vartheta_k u(t))(x)\,\D x+\kappa\int_{\Omega}\vartheta_k(x)|u(t)(x)|^2\,\D x\\&\le
(-\A_{(1-\alpha)}u(t)).(\vartheta_k u(t))+\kappa\int_{\Omega}\vartheta_k(x)|u(t)(x)|^2\,\D x+\int_\Omega  \vartheta_k (x)c(x)\,\D x,
\\&\le
\int_\Omega \bigl( (a_1C_{\overline \vartheta}^2k^{-2})|u(t)|^2+(2a_1C_{\overline\vartheta}k^{-1})|u(t)||\nabla u(t)|+(a_1C_{\vartheta}k^{-1})|u(t)||\nabla u(t)|\bigr)\D x\\& +\int_\Omega  \vartheta_k (x)c(x)\,\D x\endaligned$$
Thus, whenever $R\in \ro0,\infty..$, $\tau\in \oi0,\infty..$ and
$|u(t)|_{H^1_0}\le R$ for all $t\in\ci0,\tau..$ then
$$(V_k\circ u)'(t)+2\kappa (V_k\circ u)(t)\le b_k R^2 +c_k, \quad t\in \lo 0,\tau...$$

This implies that
$$V_k(u(t))\le e^{-2\kappa t}V_k(u(0))+(b_k R^2 +c_k)/(2\kappa),\quad t\in\ci0,\tau...$$
This clearly implies that~\rf{010207-0925} holds for every $t\in\ci0,\tau..$. This proves the theorem for $\overline\rho\in\oi2,4..$. Similar, but simpler arguments prove the theorem for $\overline\rho=2$.
The theorem is proved.\eproof
We can now prove
\proclaim{Theorem~\dft{010207-0932}}
The semiflow $\pi$ is asymptotically compact.
\endproclaim
\proof
Let $B$ be an ultimately bounded subset of $Y=H^1_0(\Omega)$, $(v_n)_n$ be a sequence in $B$ and $(t_n)_n$ be a sequence in $\ro0,\infty..$ with $t_n\to \infty$ as $n\to\infty$. We must show that there is a subsequence of $(v_n\pi t_n)_n$ which converges in $Y$.

There is a $t_B\in\ro0,\infty..$ and an $R\in\ro0,\infty..$ such that $|v\pi t|_{H^1_0}\le R$ for all $v\in B$ and $t\in\ro t_B,\infty..$. We may assume w.l.o.g. that $t_n\ge t_B+1$ for all $n\in\N$. For $n\in\N$ let $s_n=t_n-t_B$ and $u_n\co \ro0,\infty..\to H^1_0(\Omega)$ be defined by $u_n(s)=v_n\pi(t_B+s)$ for $s\in \ro0,\infty..$.
Then, for $n\in\N$, $\tau_n:=s_n-1\ge 0$ and $u_n$ is a solution of $\pi$ with $|u_n(s)|_{H^1_0}\le R$ for all $s\in\ro0,\infty..$ and $u_n(s_n)=v_n\pi t_n$.

Suppose that $\overline\rho>2$.

We claim that
$$\mycenter{There is a strictly increasing sequence $(n_m)_m$ in $\N$ and a $v\in H^1_0(\Omega)$ such that $(u_{n_m}(\tau_{n_m}))_m$ converges to $v$ in $L^2(\Omega)$.}\leqno\dff{010207-0952}$$
Let $\widetilde u_n=\phi_{(1/2),-\alpha}\circ u_n$, $n\in\N$. Then $\widetilde u_n$ is a solution of $\widetilde \pi$ and
$(\widetilde u_{n_m}(\tau_{n_m}))_m$ converges to $\widetilde v:=\phi_{(1/2),-\alpha}v$ in $\phi_{0,-\alpha}[X]$. Thus
$(\widetilde u_{n_m}(\tau_{n_m}))_m$ converges to $\widetilde v$ in $X^{-\alpha}$. Using Proposition~\rf{150107-1446} (with appropriately modified notation) we see that $(\widetilde u_{n_m}(\tau_{n_m}+1))_m$ converges to $\widetilde v\widetilde \pi 1$ in $\phi_{1/2,-\alpha}[X^{1/2}]$. This means that $( u_{n_m}(s_{n_m}))_m$ converges to $ v \pi 1$ in $X^{1/2}$ and completes the proof of the theorem.

Thus we only have to prove~\rf{010207-0952}. Let $\beta_{\text{K}}$ be the Kuratowski measure of noncompactness on $X=L^2(\Omega)$. Then for every $k\in\N$ and $n_0\in\N$,
$$\aligned\beta_{\text{K}}\{\,u_n(\tau_n)\mid n\in\N\,\}&\le
\beta_{\text{K}}\{\,(1-\vartheta_k)u_n(\tau_n)\mid n\in\N\,\}+\beta_{\text{K}}\{\,\vartheta_ku_n(\tau_n)\mid n\in\N\,\}\\&=\beta_{\text{K}}\{\,(1-\vartheta_k)u_n(\tau_n)\mid n\in\N\,\}+\beta_{\text{K}}\{\,\vartheta_ku_n(\tau_n)\mid n\ge n_0\,\}\endaligned$$
By Theorem~\rf{270606-1543} and the fact that $\tau_n\to\infty$ as $n\to\infty$, for every $\varepsilon\in\oi0,\infty..$ there are a $k\in\N$ and an $n_0\in\N$ such that $|\vartheta_k u_n|_{L^2}<\varepsilon$ for all $n\ge n_0$. Thus, for every $\varepsilon\in\oi0,\infty..$, there is a $k\in\N$ such that
$$\beta_{\text{K}}\{\,u_n(\tau_n)\mid n\in\N\,\}\le
\beta_{\text{K}}\{\,(1-\vartheta_k)u_n(\tau_n)\mid n\in\N\,\}+\varepsilon.$$
Since $1-\vartheta_k\in C^1_0(\R^N)$, results in~\cite{\rfa{PR1}} imply that the map $u\mapsto (1-\vartheta_k)u$ is compact from $H^1_0(\Omega)$ to $L^2(\Omega)$ and so
$$\beta_{\text{K}}\{\,(1-\vartheta_k)u_n(\tau_n)\mid n\in\N\,\}=0.$$ We therefore obtain that
$$\beta_{\text{K}}\{\,u_n(\tau_n)\mid n\in\N\,\}=0$$ and so there
is a strictly increasing sequence $(n_m)_m$ in $\N$ and a $v\in L^2(\Omega)$ such that $(u_{n_m}(\tau_{n_m}))_m$ converges to $v$ in $L^2(\Omega)$. Since $(u_{n_m}(\tau_{n_m}))_m$ is bounded in $H^1_0(\Omega)$, by taking a further subsequence if necessary, we may assume that $(u_{n_m}(\tau_{n_m}))_m$ converges weakly in $H^1_0(\Omega)$ (and hence in $L^2(\Omega)$) to some $w\in H^1_0(\Omega)$. Thus $w=v$ and~\rf{010207-0952} follows. This proves the theorem for $\overline\rho\in\oi2,4..$. Similar (and simpler) arguments establish the proof for $\overline\rho=2$. The proof is complete.
\eproof
We may now prove the main result of this paper.

\demo{proof of Theorem~\rf{010207-0610}}
In view of Corollary~\rf{270606-1148} and Theorem~\rf{010207-0932}, Proposition~\rf{250606-1446} implies part $(1)$ of the theorem.
If $\overline\rho= 2$,  then, noting that $\alpha=0$ in this case, Theorem~\rf{291206-1641} implies part $(2)$ of the theorem. If $\overline\rho>2$, then, by part $(1)$, $\widetilde\pi$ is a global semiflow on $\phi_{(1/2),-\alpha}[X^{1/2}]$ and it has a global attractor $\widetilde{\Cal A}:=\phi_{(1/2),-\alpha}[\Cal A]$. Thus, by Theorem~\rf{291206-1641}, $\widetilde{\Cal A}$ lies in $\phi_{1-\alpha,-\alpha}[X^{1-\alpha}]$ and is compact in $\phi_{1-\alpha,-\alpha}[X^{1-\alpha}]$. This implies part $(2)$ of the theorem in this case. The proof is complete.
\eproof
\head \dfs{Appendix} Appendix \endhead
We will prove Proposition~\rf{280107-2215}. We require the
 following simple and known result.
\proclaim{Lemma~\dft{220107-0748}}
Let $E_k$ and $F_k$, $k\in\{1,2\}$, be real or complex normed spaces with $E_2$ and $F_2$ complete. Let $e\co E_1\to E_2$, $f\co F_1\to F_2$ and $B_1\co E_1\to F_1$ be linear isometries with $B_1$ bijective. If $e[E_1]$ is dense in $E_2$ and $f[F_1]$ is dense in $F_2$, then there is a unique continuous map $B_2\co E_2\to F_2$ such that $B_2\circ e=f\circ B_1$. $B_2$ is a linear bijective isometry.
\endproclaim
\demo{Proof of Proposition~\rf{280107-2215}}
It is immediate  that for all $\alpha$, $\beta\in \R$ with $\beta\ge \alpha$ the map $\phi_{\beta,\alpha}\co X^\beta\to X^\alpha$ is defined, linear and bounded. The density of $X^\delta$ in $X^\gamma$ for all $\gamma$, $\delta\in\oi0,\infty..$ with $\delta>\gamma$ implies that
$$\mycenter {for all $\alpha$, $\beta\in \R$ with $\beta\ge \alpha$ the map $\phi_{\beta,\alpha}$ is injective and  $\phi_{\beta,\alpha}[X^\beta]$ is dense in $X^\alpha$.}\leqno\dff{210107-2130}$$
Now formula~\rf{220107-1030} and an  integration using H\"older inequality shows that
$$\mycenter{$|A^{\beta}x|_X\le |A^\alpha x|_X^{1-\theta}|A^{\gamma}x|_X^\theta$, for all $\alpha$, $\gamma\in\R$, $\delta\in\ro0,\infty..$, $\theta\in\ci0,1..$ and $x\in X^\delta$ with $\alpha\le\gamma\le\delta$ and $\beta=(1-\theta)\alpha+\theta \gamma$.}\leqno\dff{220107-1120}$$
Formula~\rf{210107-2125} and the definition of the maps $\phi_{\beta,\alpha}$ implies that
$$  \mycenter{For all $\alpha\in\oi0,\infty..$, $\beta\in\ro0,\infty..$ and $x\in X$
\par\noindent $|\phi_{\beta,\beta-\alpha }A^{-\beta}x|_{X^{\beta-\alpha}}=|A^{-\alpha}x|_X$. }\leqno\dff{220107-0825}$$
Using this formula (for $\beta=0$) we see that
$$\mycenter{$|\phi_{\delta,\alpha}x|_{X^\alpha}=|A^\alpha x|_X$, for all $\alpha\in\R$, $\delta\in\ro0,\infty..$ and $x\in X^\delta$ with $\delta\ge\alpha$.}\leqno\dff{220107-1405}$$
Now~\rf{210107-2130},~\rf{220107-1120} and~\rf{220107-1405} imply
the interpolation inequality
$$\mycenter{$|\phi_{\gamma,\beta}x|_{X^\beta}\le |\phi_{\gamma,\alpha} x|_{X^\alpha}^{1-\theta}|x|_{X^\gamma}^\theta$, for all $\alpha$, $\gamma\in\R$,  $\theta\in\ci0,1..$ and $x\in X^\gamma$ with $\alpha\le \gamma$ and $\beta=(1-\theta)\alpha+\theta \gamma$.}\leqno\dff{220107-1414}$$
A straightforward proof by cases also shows that
$$ \phi_{\alpha,\alpha}=\Id_{X^\alpha},\quad  \alpha\in\R\leqno\dff{f:frac10}$$
and
$$ \phi_{\gamma,\alpha}=\phi_{\beta,\alpha}\circ\phi_{\gamma,\beta}
, \quad\alpha,\beta,\gamma\in\R,\,
\gamma\ge\beta\ge\alpha.\leqno\dff{f:frac11}$$

For all $\alpha$, $\beta\in \ro0,\infty..$, \rf{210107-2125} implies that
 $A^{-\beta}[X^\alpha]=X^{\beta+\alpha}$ and for all $y\in X^\alpha$
 $$|A^{-\beta}y|_{X^{\beta+\alpha}}=|A^{-(\beta+\alpha)}A^\alpha y|_{X^{\beta+\alpha}}=|A^\alpha y|_X=|y|_{X^\alpha}.$$ This implies that
 $$\mycenter {For all $\alpha$, $\beta\in \ro0,\infty..$,
 $A_{(\alpha)}^{-\beta}:=A^{-\beta}|X^\alpha\co X^\alpha\to X^{\beta+\alpha}$  is a linear bijective isometry.}\leqno\dff{251006-0829}$$

For every $\alpha\in \oi0,\infty..$,
$\beta\in
\ro0,\infty..$ formula~\rf{210107-2125} implies that $\phi_{\beta,\beta-\alpha}$ is an isometry from the space $X^{\beta}$ endowed with the (in general incomplete) norm
$x\mapsto n_{\beta-\alpha}(x):= |A^{\beta-\alpha}x|_X$ to the Hilbert space $X^{\beta-\alpha}$. The same formula with $\beta=0$ shows that $\phi_{0,-\alpha}$ is an isometry from the space $X=X^{0}$ endowed with the (in general incomplete) norm
$x\mapsto n_{-\alpha}(x)= |A^{-\alpha}x|_X$ to the Hilbert space $X^{-\alpha}$. Since $A^{-\beta}$ is a bijective isometry from
$X$ endowed with the norm $n_{-\alpha}$ to $X^{\beta-\alpha}$ endowed with the norm $n_{\beta-\alpha}$ it follows from Lemma~\rf{220107-0748} that
$$\mycenter{For every $\alpha\in \oi0,\infty..$,
$\beta\in
\ro0,\infty..$  there is a unique continuous map
$A_{(-\alpha)}^{-\beta}\co X^{-\alpha}\to
X^{\beta-\alpha}$ with $A_{(-\alpha)}^{-\beta}\circ
\phi_{0,-\alpha}=\phi_{\beta,\beta-\alpha}\circ
A^{-\beta}$. $A_{(-\alpha)}^{-\beta}$ is a linear bijective isometry.}\leqno\dff{f:frac16}$$
Now, for $\alpha\in\R$ and $\beta\in\oi0,\infty..$ we may define the map
$A_{(\alpha)}^\beta\co X^\alpha\to X^{-\beta+\alpha}$ by
$$\mycenter{$A_{(\alpha)}^\beta=(A_{(-\beta+\alpha)}^{-\beta})^{-1}$. }\leqno\dff{251006-0857} $$
We also set $A_{(\alpha)}:=A_{(\alpha)}^1$.
The above definitions and simple density arguments show that
$$\mycenter{For all $\gamma$, $\gamma'\in\R$ with $\gamma>\gamma'$ and
all $\beta\in\R$,
\par\noindent
$\phi_{-\beta+\gamma,-\beta+\gamma'}\circ
A_{(\gamma)}^{\beta}=A_{(\gamma')}^{\beta}\circ
\phi_{\gamma,\gamma'}$
}\leqno\dff{251006-0858} $$
and
$$\mycenter{For all $\alpha$, $\beta$ and $\gamma\in\R$,
$A_{(-\gamma+\alpha)}^\beta\circ A_{(\alpha)}^\gamma=A_{(\alpha)}^{\beta+\gamma}$.}\leqno\dff{251006-0926} $$
Since for $\alpha$, $\beta\in\R$ with $\beta\ge \alpha$ the map $\phi_{\beta,\alpha}$ is bijective onto its range, we may,  for all $\alpha$, $\beta\in\oi0,\infty..$ define the map
$$\widetilde A_{(-\alpha)}^\beta:=A_{(\beta-\alpha)}^\beta\circ \phi_{\beta-\alpha,-\alpha}^{-1}\co \phi_{\beta-\alpha,-\alpha}[X^{\beta-\alpha}]\subset X^{-\alpha}\to X^{-\alpha}.$$ We also set $\widetilde A_{(-\alpha)}:=\widetilde A_{(-\alpha)}^1$.

It follows that $\widetilde A_{(-\alpha)}^\beta$ is bijective and its inverse is
 $\widetilde
A_{(-\alpha)}^{-\beta}:=\phi_{\beta-\alpha,-\alpha}\circ A_{(-\alpha)}^{-\beta}$.

 We claim that
$\widetilde A_{(\beta-\alpha)}^{\beta}$ is symmetric with respect to the scalar product on $X^{-\alpha}$.
First notice that
$$R_\alpha^{-1}\phi_{0,-\alpha}x=A^{-2\alpha}x,\quad \alpha\in\oi0,\infty..,\,x\in X.\leqno\dff{220107-1912}$$
This implies that, for $u$ and $v\in X$
$$\langle \phi_{0,-\alpha}u,\phi_{0,-\alpha}v\rangle_{X^{-\alpha}}=\langle A^{-2\alpha}u, A^{-2\alpha}v\rangle_{X^\alpha}=\langle A^{-\alpha}u, A^{-\alpha}v\rangle_{X} $$ so
$$\langle \phi_{0,-\alpha}u,\phi_{0,-\alpha}v\rangle_{X^{-\alpha}}=\langle A^{-\alpha}u, A^{-\alpha}v\rangle_{X}, \quad \alpha\in\oi0,\infty..,\,u,\,v\in X.\leqno\dff{220107-1917}$$
Since $ A_{(\beta-\alpha)}^{\beta}$ is the inverse of
$ A_{(-\alpha)}^{-\beta}$ and $A^\beta$ is the inverse of $A^{-\beta}$, we obtain from~\rf{f:frac16}
$$A_{(\beta-\alpha)}^\beta\circ \phi_{\beta,\beta-\alpha}=\phi_{0,-\alpha}\circ A^\beta.\leqno\dff{230107-0914}$$
Hence, for $u\in X^\beta$, we obtain from~\rf{230107-0914},
$$\widetilde A_{(\beta-\alpha)}^{\beta}\phi_{\beta-\alpha,-\alpha}\phi _{\beta,\beta-\alpha}u=A_{(\beta-\alpha)}^{\beta}\phi _{\beta,\beta-\alpha}u
=\phi_{0,-\alpha}A^{\beta}u.$$ Using~\rf{220107-1917}, we thus obtain, for $u$, $v\in X^\beta$,  $x=\phi_{\beta-\alpha,-\alpha}\phi _{\beta,\beta-\alpha}u$ and $y=\phi_{\beta-\alpha,-\alpha}\phi _{\beta,\beta-\alpha}v$
$$\langle \widetilde A_{(\beta-\alpha)}^{\beta}x, y\rangle_{X^{-\alpha}}=\langle A^{-\alpha}u, A^{-\alpha}A^{\beta}v\rangle_X=\langle A^{-\alpha}u, A^{\beta}A^{-\alpha}v\rangle_X $$ and in particular
$$\mycenter{For $u\in X^\beta$ and $x=\phi_{\beta-\alpha,-\alpha}\phi _{\beta,\beta-\alpha}u$,\par\noindent\hfill$\langle \widetilde A_{(\beta-\alpha)}^{\beta}x, x \rangle_{X^{-\alpha}}=\langle A^{-\alpha}u, A^{\beta}A^{-\alpha}u\rangle_X$.\hfill}\leqno\dff{230107-0938} $$
Similarly,
$$\langle x,\widetilde A_{(\beta-\alpha)}^{\beta}y \rangle_{X^{-\alpha}}=\langle  A^{\beta}A^{-\alpha}u,A^{-\alpha}v\rangle_X. $$
Now the symmetry of $A^{\beta}$ relative to the scalar product on $X$  shows that
$$\langle \widetilde A_{(\beta-\alpha)}^{\beta}x, y\rangle_{X^{-\alpha}}=\langle x, \widetilde A_{(\beta-\alpha)}^{\beta}y \rangle_{X^{-\alpha}}.\leqno\dff{220107-1936}$$
for $x$, $y\in \phi_{\beta-\alpha,-\alpha}[\phi_{\beta,\beta-\alpha}[X^\beta]]$. Thus, by density,~\rf{220107-1936} holds for all $x$, $y\in \phi_{\beta-\alpha,-\alpha}[X^{\beta}]$, proving our claim. This claim implies that the map $\widetilde A_{(\beta-\alpha)}^{\beta}$ is self-adjoint and so its inverse $\widetilde A_{(-\alpha)}^{-\beta}$ is symmetric on $X^{-\alpha}$.
Moreover,~\rf{251006-0926} implies that
$$\widetilde A_{(-\alpha)}^{-\beta-\gamma}=\widetilde A_{(-\alpha)}^{-\beta}\circ \widetilde A_{(-\alpha)}^{-\gamma},\quad \alpha,\,\beta,\,\gamma\in\oi0,\infty...\leqno\dff{230107-0953}$$
In particular, $\widetilde A_{(-\alpha)}^{-\beta}$ is nonnegative.

Now, let $\alpha\in\oi0,\infty..$ be arbitrary. Since $\re\sigma(A)>0$ there is a $\delta\in\oi0,\infty..$ such that
$$\langle Ax,x\rangle_X\ge \delta\langle x,x\rangle_X,\quad x\in X.$$
Hence, by~\rf{230107-0938}, for $u\in X^1$ and $x=\phi_{1-\alpha,-\alpha}\phi _{1,1-\alpha}u$,
$$\langle \widetilde A_{(1-\alpha)}^{1}x, x \rangle_{X^{-\alpha}}=\langle A^{-\alpha}u, A^{1}A^{-\alpha}u\rangle_X\ge \delta\langle A^{-\alpha}u,A^{-\alpha}u\rangle_X=\delta\langle x,x\rangle_{X^{-\alpha}} \leqno\dff{230107-1002} $$
This implies, by density, that $\re\sigma(\widetilde A_{(1-\alpha)}^{1})>0$ so $B:=\widetilde A_{(1-\alpha)}$ generates the family $B^{-\beta}$, $\beta\in\oi0,\infty..$, of basic fractional power spaces of $B$. By~\rf{210107-2125}
$$B^{-\beta-\gamma}=B^{-\beta}\circ B^{-\gamma},\quad \beta,\,\gamma\in\oi0,\infty...\leqno\dff{230107-1048}$$ We claim that
$$B^{-\beta}=\widetilde A_{(-\alpha)}^{-\beta},\quad \beta\in\oi0,\infty...\leqno\dff{230107-1050}$$
Let $Z$ be the set of $\beta\in\oi0,\infty..$ with $B^{-\beta}=\widetilde A_{(-\alpha)}^{-\beta}$. Clearly, $1\in Z$ and so induction on $m\in \N$ using~\rf{230107-0953} and~\rf{230107-1048} imply that $Z$ contains all integers. Since nonnegative symmetric operators on a Hilbert space have unique nonnegative square roots, it follows by induction on $k\in\N$ (again using \rf{230107-0953} and~\rf{230107-1048}) that
$Z$ contains all numbers of the form $m/2^k$ with $m$, $k\in\N$. The set $Z_0$ of such numbers is dense in $\oi0,\infty..$. Now let $\beta\in\oi0,\infty..$ be arbitrary and $(\beta_n)_n$ be a sequence in $Z_0$ converging to $\beta$.
By formula~\rf{230107-2127}  we have that
$$|B^{-\beta_n}x- B^{-\beta}x|_{X^{-\alpha}}\to 0,\quad x\in X^{-\alpha}$$
and
$$|A^{-\beta_n}x- A^{-\beta}x|_{X}\to 0,\quad x\in X.$$
In particular, using the fact that
$$\widetilde A_{-\alpha}^{-\gamma}\phi_{0,-\alpha}u=\phi_{0,-\alpha}A^{-\gamma}u,\quad \alpha,\,\gamma\in\oi0,\infty..,\,u\in X$$ we obtain that
$$|B^{-\beta_n}\phi_{0,-\alpha}u- B^{-\beta}\phi_{0,-\alpha}u|_{X^{-\alpha}}\to 0,\quad u\in X$$
$$|\widetilde A_{-\alpha}^{-\beta_n}\phi_{0,-\alpha}u- \widetilde A_{-\alpha}^{-\beta}\phi_{0,-\alpha}u|_{X^{-\alpha}}\to 0,\quad u\in X$$
Thus
$$B^{-\beta}\phi_{0,-\alpha}u=\widetilde A_{-\alpha}^{-\beta}\phi_{0,-\alpha}u,\quad u\in X$$
so, by density, $B^{-\beta}=\widetilde A_{-\alpha}^{-\beta}$ and thus $\beta\in Z$. This proves our~\rf{230107-1050}.
We obtain from~\rf{230107-1050} that
$$X^\beta_B=\phi_{\beta-\alpha,-\alpha}[X^{\beta-\alpha}].$$
We also claim that
$$\mycenter{For all $\beta\in\oi0,\infty..$, $\phi_{\beta-\alpha,-\alpha}$ is an isometry of $X^{\beta-\alpha}$ onto $X^\beta_B$.}\leqno\dff{280107-1610}$$
To prove this claim, let $\widetilde x$, $\widetilde y\in X^{\beta-\alpha}$ be arbitrary and let $x=\phi_{\beta-\alpha,-\alpha}\widetilde x$, $y=\phi_{\beta-\alpha,-\alpha}\widetilde y$. Suppose first that $\widetilde x=\phi_{\beta,\beta-\alpha}u$, $\widetilde y=\phi_{\beta,\beta-\alpha}v$ with $u$, $v\in X^\beta$. Then, by~\rf{230107-0914},
$B^\beta x=\phi_{0,-\alpha}A^\beta u$ and $B^\beta y=\phi_{0,-\alpha}A^\beta v$. Therefore, using~\rf{220107-1917} we obtain
$$\aligned\langle x,y\rangle_{X^\beta_B}&=\langle B^\beta x,B^\beta y\rangle_{X^{-\alpha}}=\langle\phi_{0,-\alpha} A^\beta u,\phi_{0,-\alpha}A^\beta v\rangle_{X^{-\alpha}}=\langle A^{-\alpha}A^\beta u,A^{-\alpha}A^\beta v\rangle_X\\&=\langle A^{\beta-\alpha}u,A^{\beta-\alpha}v\rangle_X.\endaligned$$
If $\beta-\alpha\ge 0$, then $$\langle A^{\beta-\alpha}u,A^{\beta-\alpha}v\rangle_X=\langle u,v\rangle_{X^{\beta-\alpha}}=\langle \widetilde x,\widetilde y\rangle_{X^{\beta-\alpha}}.$$
If $\beta-\alpha<0$, then, by~\rf{220107-1917},
$$\aligned\langle A^{\beta-\alpha}u,A^{\beta-\alpha}v\rangle_X&=\langle A^{-(\alpha-\beta)}u,A^{-(\alpha-\beta)}v\rangle_X=\langle \phi_{0,-(\alpha-\beta)}u,\phi_{0,-(\alpha-\beta)}v\rangle_{X^{-(\alpha-\beta)}}\\&=
\langle \widetilde x,\widetilde y\rangle_{X^{-(\alpha-\beta)}}\endaligned$$
Thus, in both cases,
$\langle x,y\rangle_{X^\beta_B}=\langle \widetilde x,\widetilde y\rangle_{X^{\beta-\alpha}}$.
Since the set $\phi_{\beta,\beta-\alpha}[X^\beta]$ is dense in $X^{\beta-\alpha}$,~\rf{280107-1610} follows.

We further claim that
$$\mycenter{Whenever $\alpha\in\ro0,(1/2)..$, $x\in X^{1-\alpha}$ and $v\in X^{1/2}\subset X^\alpha$, then \par\hfill$(A_{(1-\alpha)}x).v=\langle  x,v\rangle_{X^{1/2}}$.\hfill}\leqno\dff{280107-1659}$$
Here, the dot `$.$' denotes function application between an element of $X^{-\alpha}$ and $X^\alpha$.
To prove this claim, assume first that $ x\in X^1$. Then, by~\rf{230107-0914}, $ A_{(1-\alpha)}x=\phi_{0,-\alpha}A x$ and so using the Fr\'echet-Riesz theorem and~\rf{220107-1912} we obtain
$$(A_{(1-\alpha)}x).v=\langle v,R_\alpha^{-1}\phi_{0,-\alpha}A x\rangle_{X^\alpha}=\langle v,A^{-2\alpha}A x\rangle_{X^\alpha}=\langle v,A x\rangle _X=\langle v, x\rangle_{X^{1/2}}$$
so the claim follows in this special case. The general case follows by the density of $X^1$ in $X^{1-\alpha}$ and in $X^{1/2}$.

The proposition is proved.
\eproof

\enddocument